\documentclass{amsproc}
\usepackage{amssymb}

\usepackage{graphicx}
\usepackage{amscd}

\theoremstyle{plain}

\newtheorem{corollary}{Corollary}

\newtheorem{lemma}{Lemma}

\newtheorem{proposition}{Proposition}
\newtheorem{remark}{Remark}

\newtheorem{theorem}{Theorem}
\numberwithin{equation}{section}

\input{tcilatex}
 
\begin{document}
\title[HIGGS]{HIGGS BUNDLES AND HOLOMORPHIC FORMS}
\author{Walter Seaman}
\address{Department of Mathematics University of Iowa}
\email{walter-seaman@uiowa.edu}
\urladdr{http://www.math.uiowa.edu/\symbol{126}seaman}
\date{spring 1999}

\begin{abstract}
For a complex manifold $X$ which has a holomorphic form $\varpi $ of odd
degree $k$, we endow $E^{a}=\bigoplus_{p\geq a}\Lambda ^{(p,0)}(X)$ with a
Higgs bundle structure $\theta $\ given by $\theta (Z)(\phi ):=\{i(Z)\varpi
\}\wedge \phi $. The properties such as curvature and stability of these and
other Higgs bundles are examined. We prove (Theorem 2, section 2, for $k>1$) 
$E^{a}$ and additional classes of Higgs subbundles of $E^{a}$ do not admit
Higgs-Hermitian-Yang-Mills metric in any one of the cases: i. $\deg (X)<0$,
ii. $\deg (X)=0$ and $a\leq n-k+1$, or iii. $a\leq n-k+1$ and $k\geqslant 
\frac{n}{2}+1$. \ We give examples of (noncompact) K\"{a}hler manifolds with
the above Higgs structure which admit Higgs-Hermitian-Yang-Mills metrics. We
also examine vanishing theorems for $(p,q)-$forms with values in Higgs
bundles.
\end{abstract}

\maketitle

\section{Section 1}

\refstepcounter{section}

The purpose of this paper is to give new examples of Higgs bundles which
arise in a rather natural way, and to study their properties. Recall that a $%
Higgs$ $Bundle$ \cite{simpIHES} is a holomorphic vector bundle, $E$ $%
\longrightarrow X$ over a complex manifold $X$, together with a holomorphic
section $\theta \in \vartheta \Gamma (Hom(E)\otimes \Lambda ^{1,0}(X))$
,(the ``Higgs'' form) which satisfies the equation $\theta \wedge \theta =0.$
This equation means that if $Z$ and $W$ are holomorphic tangent vectors to X
at a point, then $[\theta (Z),\theta (W)]$ = 0 as an endomorphism of $E$ \
at that point.

The examples consist of a complex manifold $X$ of complex dimension n which
is assumed to possess a nontrivial holomorphic k-form $\varpi $ where $k$ is
odd. The bundle $E$ is given by $E:=\overset{n}{\underset{p=0}{\bigoplus }}%
\bigwedge^{(p,0)}(X)$ , and the Higgs form $\theta $ is given by the
prescription $\theta (Z)(\phi ):=\{i(Z)\varpi \}\wedge \phi $, where $\phi $
is a section of $E$ and $Z$ is a holomorphic tangent vector. Defining $E^{a} 
$ by $E^{a}:=\overset{n}{\underset{p=a}{\bigoplus }}\bigwedge^{(p,0)}(X)$ ( $%
E=$ $E^{0}$), the $E^{a}$ form a Higgs filtration of $E$ (cf. \ref{filt1}).
We now give some examples of complex manifolds possessing such forms.

i. $X=$ any complex torus.

ii. If $X$ is the zero-locus in $\mathbb{P}^{n+1}$of a homogeneous
polynomial of large degree D, then $h^{n,0}(X)=\binom{D-1}{n+1}$ so if n is
odd these are examples of the types of complex manifolds required.

iii. Calabi-Yau manifolds-compact K\"{a}hler Ricci flat complex 3-manifolds
with a nowhere vanishing holomorphic 3-form, i.e. trivial canonical bundle,
and higher-dimensional analogs (cf. \cite{bryant}, pages 144-145).

iv. For any complex manifold $X$, its holomorphic cotangent bundle $\Lambda
^{(1,0)}X$ admits a canonical holomorphic one-form $\phi $ $\in \vartheta
\Gamma $ $\Lambda ^{(1,0)}(\Lambda ^{(1,0)}X)$ such that $\partial \phi $ is
a (holomorphic) symplectic two-form. This $\phi $ can be given invariantly
by the formula $\phi (Z_{\alpha })=\alpha (\pi _{\ast \alpha }(Z_{\alpha
})),Z_{\alpha }\in T_{\alpha }^{(1,0)}(\Lambda ^{(1,0)}X),\alpha \in \Lambda
^{(1,0)}X$ with $\pi :\Lambda ^{(1,0)}X\rightarrow X$ the projection (cf. 
\cite{bryant}, pages 85-86). \ Replacing $X$ with the complex manifold $%
\Lambda ^{(1,0)}X$, one gets the corresponding holomorphic one-form $\Phi $ $%
\in \vartheta \Gamma $ $(\Lambda ^{(1,0)}(\Lambda ^{(1,0)}(\Lambda
^{(1,0)}X))$ and (symplectic) two-form $\partial \Phi \in \vartheta \Gamma $ 
$(\Lambda ^{(2,0)}(\Lambda ^{(1,0)}(\Lambda ^{(1,0)}X))$ on $\Lambda
^{(1,0)}(\Lambda ^{(1,0)}X)$. Let $p:$ $\Lambda ^{(1,0)}(\Lambda
^{(1,0)}X)\rightarrow \Lambda ^{(1,0)}X\ $be the projection. \ Then for any
holomorphic functions $a$ and $b$ on $\Lambda ^{(1,0)}(\Lambda ^{(1,0)}X)$,
one gets a holomorphic three-form $a\Phi \wedge p^{\ast }\partial \phi +$ $%
b\partial \Phi \wedge p^{\ast }\phi \in \vartheta \Gamma $ $(\Lambda
^{(3,0)}(\Lambda ^{(1,0)}(\Lambda ^{(1,0)}X))$). Computation of these
3-forms in local holomorphic coordinates (using coordinates on $\Lambda
^{(1,0)}X$ given by ''pulling up'' a holomorphic chart on $X$ and then
''pulling up'' these coordinates on $\Lambda ^{(1,0)}X$ via $p$ to $\Lambda
^{(1,0)}(\Lambda ^{(1,0)}X)$) shows that these forms are generally nonzero.

v. If $M$ is any of the above examples, then any complex manifold $%
\widetilde{M\text{ }}$from which there is a holomorphic submersion $p:%
\widetilde{M\text{ }}\rightarrow M$ onto $M$, itself inherits nonzero
holomorphic odd-degree forms from $M$ by pull-back. \ For example, coverings
or blowing up any of the above examples at any number of points and/or
taking products of those examples will serve as such an $\widetilde{M\text{ }%
}$.

We investigate the curvature, stability and other properties of these Higgs
bundles (and also general Higgs Bundles) and prove the following :

\bigskip

\textbf{Theorem} \ref{bigthm}, section 2. Let $X$ be a compact K\"{a}hler
manifold with a nontrivial holomorphic k-form $\varpi $ where $k>1$ is odd.
Let the Higgs structure of $E$ be as above. and let $P$ be any Higgs
subbundle of $E$ of the form $P=\bigoplus_{s=1}^{z}\Lambda ^{(p_{s},0)}(X)$, 
$0\leq p_{1}<p_{2}<\cdots <p_{z}\leq n$, $(z\geq 2)$. Then $P$ does not
admit any Higgs-Hermitian-Yang-Mills metric in any of the following cases :

i. $\deg (X)<0$

ii. $\deg (X)=0$ and $p_{1}\leq n-k+1$

iii. $k\geq \frac{n}{2}+1$, $p_{1}\leq n-k+1$, and $\varpi $ is a section of 
$P$.\bigskip

Note that the degree statement in ii. is sharp because the Higgs form $%
\theta $ acts trivially on $E^{n-k+2}$. If $X$ is compact K\"{a}hler with
first Chern class $c_{1}(X)=0$, then the Yau resolution of the Calabi
conjecture \cite{Yau}, yields a Ricci-flat metric $g\ $on $X.$ Extending $g$
in the usual way to the complex exterior algebra of $X$ gives a Hermitian
metric on $E^{n-k+2}$ which is Higgs-Hermitian-Yang-Mills in the ''vacuous''
sense that $g$ is Hermitian-Yang-Mills and the Higgs form vanishes.

We also examine Bochner-type vanishing results (\ref{thmw=0}, section 2) and
Kodaira--Nakano-type vanishing theorems (\ref{thm1} and \ref{thm3.2}) in
this setting.

The original study of Higgs bundles is due to Hitchin \cite{hitchin}, where
the case of rank 2 vector bundles over curves is considered. Hitchin studies
the Yang-Mills equations with ``interaction term'' given by the Higgs field
(cf. the discussion above \ref{c1}). Hitchin obtains a correspondence
relating irreducible rank 2 flat vector bundles and degree zero stable Higgs
bundles over Riemann surfaces. This correspondence has its genesis in the
work of Narasimhan and Seshardi \cite{narasihiman}.

Higgs bundles also arise in the study of \ Variations of Hodge Structure.
See e.g. \cite{schmidt} sections 1 and 2, \cite{green},\cite{wells} Chapter
V section 6, \cite{simpConstVHS}pages 868-869, and \cite{simpUbiquity}
section 1, for detailed information. \ Generalizing the idea that Hitchin
had introduced, Simpson \cite{simpConstVHS},\cite{simpUbiquity}, \cite
{simpIHES} defined the notion of Higgs bundles on higher dimensional
varieties, where the equation $\theta \wedge \theta =0$ (automatically
satisfied on a curve) is part of the definition. Simpson studied the moduli
space of stable Higgs bundles with vanishing Chern classes in work which
leads up to the following striking result (showing the ``ubiquity'' of VHS)
among others:

\textit{If }$M$\textit{\ is a smooth projective variety then any
representation of }$\pi _{1}(M)$\textit{\ can be deformed to a
representation arising from a complex variation of Hodge structure}.\newline
This result, among other things, restricts the types of groups which can
arise as the fundamental group for any such $M,$ cf. \cite{amoros} chapter 7.

The author would like to thank Aroldo Kaplan, Carlos Simpson and Olivier
Debarre for their help in various aspects of this work.

We now continue with the development of properties of Higgs bundles. Any
Higgs bundle has a naturally defined operator $D^{^{\prime \prime }}:\Gamma
E\rightarrow \Gamma (E\otimes \Lambda ^{1}(X))$defined by $D^{\prime \prime
}=\overset{\_}{\partial }+\theta $ where $\overset{\_}{\partial }$ is the
complex structure on $E$ .The three conditions: $\overset{\_}{\partial }$ is
integrable ($\overline{\partial }^{2}=0),$ $\theta $ is holomorphic and $%
\theta \wedge \theta =0$ are simultaneously expressed in the single equation 
$(D^{^{\prime \prime }})^{2}=0$.

Let $h$ be a Hermitian metric on $E$. The Hermitian connection of $(E,h)$ ,$%
\nabla ,$ can be uniquely written $\nabla =\partial _{h}+\overset{\_}{%
\partial }$. Define the Hermitian adjoint of $\theta ,\overline{\theta _{h}}$
by the formula

\begin{equation}
h(\overline{\theta _{h}}(Y)\ s,t)=h(s,\theta (\ \overline{Y}\ )\ t)
\end{equation}
where $Y$ is a complex tangent vector and $s$ and $t$ are sections of $E$ .
Define $D_{h}^{^{\prime }}$ by

\begin{equation}
D_{h}^{^{\prime }}=\partial _{h}+\overline{\theta _{h}}  \label{D'h}
\end{equation}
and

\begin{equation}
D_{h}=D_{h}^{^{\prime }}+D^{^{\prime \prime }}  \label{Dh}
\end{equation}

One checks that $(D_{h}^{^{\prime }})^{2}=0$ ,that $D_{h}$ is a connection
on $E$ and that the curvature of $D_{h}$ is given by

\begin{equation}
F_{h}=(D_{h})^{2}=D_{h}^{^{\prime }}D^{^{\prime \prime }}+D^{^{\prime \prime
}}D_{h}^{^{\prime }}  \label{F}
\end{equation}

Let $\Theta =\nabla ^{2}$ be the curvature of $h.$ Although $\Theta $ is a
type (1,1) $End(E)$-valued form, in general $F_{h}$will have parts of type $%
(2,0),(1,1)$ and $(0,2).$ The relation between the components of $F_{h}$,$%
\Theta ,\theta $ will now be described. Let $\{e_{\alpha }\}_{\alpha =1}^{r}$
be a local holomorphic frame for $E$ ($r=$rank of $E)$, $h_{\alpha \overline{%
\beta }}=h(e_{\alpha },e_{\beta })$, and ($h^{\beta \overline{\gamma }})$ be
the inverse matrix of $(h_{\alpha \overline{\beta }})$. Then $\nabla
e_{\alpha }=\overset{r}{\underset{\beta =1}{\sum }}$ $e_{\beta }\bigotimes
C_{\alpha }^{\beta }$, where $C=h^{-1}\partial h$ and $\Theta =\overline{%
\partial }C.$ Also, $\theta e_{\alpha }=\overset{r}{\underset{\beta =1}{\sum 
}}e_{\beta }\bigotimes \theta _{\alpha }^{\beta }$ where $\theta _{\alpha
}^{\beta }$ are the matrix representative $(1,0)$-forms of $\theta $
relative to $\{e_{\alpha }\}.$ In this setting we also have $\overline{%
\theta _{h}}e_{\alpha }=\overset{r}{\underset{\beta =1}{\sum }}e_{\beta
}\bigotimes \overline{\theta _{h}}_{\alpha }^{\beta }$ , where $\overline{%
\theta _{h}}_{\alpha }^{\beta }=\underset{\gamma ,\kappa }{\sum }h^{\beta 
\overline{\gamma }}\;\overline{\theta _{\gamma }^{\kappa }}\;h_{\alpha 
\overline{\kappa }}$, If the frame $\{e_{\alpha }\}_{\alpha =1}^{r}$ is
orthonormal at a point of evaluation, then $\overline{\theta _{h}}_{\alpha
}^{\beta }=\overline{\theta _{\beta }^{\alpha }}$ at that point.

Now, writing $F_{h}e_{\alpha }=\underset{\beta }{\sum }\{e_{\beta }\otimes
F_{h\alpha }^{\beta }\},$ and ($F_{h}e_{\alpha })^{(a,b)}=\underset{\beta }{%
\sum }\{e_{\beta }\otimes (F_{h\alpha }^{\beta })^{(a,b)}\}$ where $%
(a,b)=(2,0),(0,2)$ or $(1,1)$ one computes (cf. \cite{simpConstVHS} page
879, fourth line from the top and also Proposition \ref{prop} below)

\begin{align}
F_{h}^{(2,0)}& =\partial \ \theta +C\wedge \theta +\theta \wedge C
\label{F20} \\
F_{h}^{(0,2)}& =\overline{\partial \ }\;\overline{\theta _{h}}  \label{F02}
\\
F_{h}^{(1,1)}& =\Theta +\theta \wedge \overline{\theta _{h}}+\overline{%
\theta _{h}}\wedge \theta  \label{F11}
\end{align}
in more detail

\begin{align*}
(F_{h}e_{\alpha })^{(2,0)}& =\underset{\beta }{\sum }\{e_{\beta }\otimes \{\
\partial \theta _{\alpha }^{\beta }+\underset{\gamma }{\sum }(C_{\gamma
}^{\beta }\wedge \theta _{\alpha }^{\gamma }+\theta _{\gamma }^{\beta
}\wedge C_{\alpha }^{\gamma })\}\} \\
(F_{h}e_{\alpha })^{(0,2)}& =\underset{\beta }{\sum }\{e_{\beta }\otimes \{%
\overline{\partial }\text{ \ }\overline{\theta _{h}}_{\alpha }^{\beta }\}\}
\\
(F_{h}e_{\alpha })^{(1,1)}& =\underset{\beta }{\sum }\{e_{\beta }\otimes
\{\Theta _{\alpha }^{\beta }+\underset{\gamma }{\sum }(\theta _{\gamma
}^{\beta }\wedge \ \overline{\theta _{h}}_{\alpha }^{\gamma }+\overline{%
\theta _{h}}_{\gamma }^{\beta }\wedge \theta _{\alpha }^{\gamma })\}\}
\end{align*}

In the course of proving \ref{F11}$\ $\ one must use the identity 
\begin{equation*}
\partial \overline{\theta _{h}}+C\wedge \overline{\theta _{h}}+\overline{%
\theta _{h}}\wedge C=0
\end{equation*}
which in turn follows from the identities $C=h^{-1}\partial h$ and $%
\overline{\theta _{h}}=h^{-1}\overline{\theta }h.$

If $M$ is any smooth manifold and $\mathbb{V\rightarrow }M$ is any real or
complex vector bundle with a connection $\nabla $: $C^{\infty }\Gamma 
\mathbb{V\rightarrow }C^{\infty }\Gamma (\mathbb{V\otimes }\Lambda ^{1}(V))$
there is a ``natural'' extension of $\nabla $ , $d^{\nabla }:\mathbb{V}%
^{k}\rightarrow \mathbb{V}^{k+1}$ , where $\mathbb{V}^{r}:=C^{\infty }\Gamma
(\mathbb{V\otimes }\Lambda ^{r}(V))$. This implies the following (cf. \cite
{simpConstVHS}, page 879):

\begin{proposition}
\label{prop}$F_{h}^{(2,0)}=d^{\nabla }\theta $ and $F_{h}^{(0,2)}=d^{\nabla }%
\overline{\theta _{h}}$.
\end{proposition}

At any point p we can always find a local holomorphic frame $\{e_{\alpha
}\}_{\alpha =1}^{r}$ adapted to p, and also $\overset{}{\underset{}{%
\overline{\partial }\text{ }\overline{\text{ }\theta _{\beta }^{\alpha }}}}=%
\overline{\overset{}{\underset{}{\partial \theta _{\beta }^{\alpha }}}}$ so
we conclude $d^{\nabla }\theta (p)=0\Longleftrightarrow $ $d^{\nabla }%
\overline{\theta _{h}}(p)=0.$ Now the above Proposition \ref{prop} implies

\begin{equation}
F_{h}^{(2,0)}=0\Longleftrightarrow F_{h}^{(0,2)}=0\Longleftrightarrow
d^{\nabla }\theta =0\Longleftrightarrow d^{\nabla }\overline{\theta _{h}}=0
\end{equation}

We now examine the curvature terms appearing in \ref{F11}. If $Z,W\ $are
holomorphic tangent vectors at a point , then \ref{F11} implies 
\begin{equation}
F_{h}(Z,\overline{W})s=\Theta (Z,\overline{W})s+\theta (Z)\overline{\theta
_{h}}(\overline{W})s-\overline{\theta _{h}}(\overline{W})\theta (Z)s=\Theta
(Z,\overline{W})s+[\theta (Z)\overline{,\theta _{h}}(\overline{W})]s
\label{crv}
\end{equation}
where s is any section of $E$. Relative to the local framing $\{e_{\alpha
}\}_{\alpha =1}^{r}$ of $E$, \ref{crv} can be written

\begin{equation}
F_{h}(Z,\overline{W})e_{\alpha }=\overset{}{\underset{}{\underset{\beta }{%
\sum }\{e_{\beta }\otimes \{\Theta _{\alpha }^{\beta }(Z,\overline{W})+%
\underset{\gamma }{\sum }}}(\theta _{\gamma }^{\beta }(Z)\ \overline{\theta
_{h}}_{\alpha }^{\gamma }(\overline{W})-\overline{\theta _{h}}_{\gamma
}^{\beta }(\overline{W})\theta _{\alpha }^{\gamma }(Z))\}\}
\end{equation}
One final identity we will use following from \ref{crv} is:

\begin{equation}
h(F_{h}(Z,\overline{Z})s,s)=h(\Theta (Z,\overline{Z})s,s)+\Vert \overline{%
\theta _{h}}(\overline{Z})s\Vert _{h}^{2}-\Vert \theta (Z)s\Vert _{h}^{2}
\label{bochner}
\end{equation}

One can see an earliest version of this formula in \cite{schmidt}, section 7
and especially Lemma (7.18), pages 271-272. In the (VHS) context of that
paper one would have $F_{h}=0.$

If we now endow $X$ with a Hermitian metric $g$, then use $g$ to take the
trace of the identity of \ref{bochner} in the $``Z"$ variables we get

\begin{equation}
h(i\Lambda F_{h}s,s)=h(i\Lambda \Theta s,s)+\overset{n}{\underset{i=1}{\sum }%
\{}\Vert \overline{\theta _{h}}(\overline{Z}_{i})s\Vert _{h}^{2}-\Vert
\theta (Z_{i})s\Vert _{h}^{2}\}  \label{boch2}
\end{equation}
where the \{$Z_{i}\}_{i=1}^{n}$ forms an orthonormal basis for $T^{1,0}X$ at
a point and also in \ref{boch2} we have used the term $i\Lambda $ as a
shorthand for ``trace with respect to $g$ over $T^{(1,0)}X$''. This can be
written\ , for example, $i\Lambda \Theta =\sum_{i=1}^{n}\Theta (Z_{i},%
\overline{Z_{i}}),$ where $\{Z_{i}\}_{i=1}^{n}$ is a $g-$orthonormal basis
for $T^{(1,0)}X$ at a point, and in general $\Lambda =i\sum g^{i\overline{j}%
}i(Z_{i})i(\overline{Z}_{j})$. If $g$ happens to be a K\"{a}hler metric then
this agrees with the usual symbolism. If $s$ is a holomorphic section of $E$%
, then we have the well-known identity (\cite{kobayashi} Chapter III
Proposition 1.5, page 50 and \cite{wu} page 349, item (3.30))

\begin{equation}
h(i\Lambda \Theta s,s)=-i\Lambda \partial \overline{\partial }\Vert s\Vert
_{h}^{2}+\Vert \nabla s\Vert _{h}^{2}  \label{holosect}
\end{equation}

Now \ref{boch2} and \ref{holosect} lead, via the Bochner technique, to the
following vanishing result (cf. \cite{wu} Theorem 5, pages 347-349, \cite
{kobayashi} Chapter III Theorem 1.9,page 52,\cite{schmidt} Lemma (7.18) ,
pages 271-272)

\begin{lemma}
\label{lemma1}Suppose $X$ is compact, $s$ is a holomorphic section of $E$
satisfying $\theta s=0$ and $i\Lambda F_{h}\leq 0$ ( pointwise \ as an
endomorphism of $E$\textit{)}. Then $s$ is parallel $\nabla s=0$ and
satisfies $\overline{\theta _{h}}s=0$ and $i\Lambda F_{h}(s)=0.$ If $%
i\Lambda F_{h}$ is a \textit{quasinegative operator (\cite{wu}, page 323)
then }$s=0$.
\end{lemma}

We will see in section 3 how Lemma \ref{lemma1} extends to
Kodaira-Nakano-type vanishing result for $(p,q)-$forms with values in a
Higgs bundle.

\section{Section 2}

\refstepcounter{section}

Let $X$ be a complex manifold of complex dimension $n$. Let $E:=\overset{n}{%
\underset{p=0}{\bigoplus }}\bigwedge^{(p,0)}(X)$ be the holomorphic vector
bundle of forms of degree $(p,0)$ for all $p$. Assume $X$ has a holomorphic
form $\varpi $ (not everywhere zero) of type $(k,0)$ where $k$ is \textit{%
odd.\ }We define a Higgs form $\theta $ on $E,\theta \in \vartheta \Gamma
(Hom(E)\otimes \Lambda ^{1,0}(X)),$ by the formula:

\begin{equation}
\theta (Z)(\phi ):=\{i(Z)\varpi \}\wedge \phi  \label{higgsform}
\end{equation}
where $Z$ is a complex tangent vector to $X$, $i(Z)$ is interior
multiplication by $Z$, and $\phi $ is any section of $E$. One can write $%
\theta $ without referring to a specific complex tangent vector locally by
the formula

\begin{equation}
\theta (\phi ):=\overset{n}{\underset{i=1}{\sum }(}\{i(\frac{\partial }{%
\partial z_{i}})\varpi \}\wedge \phi \})\otimes dz_{i}  \label{higgsform2}
\end{equation}
where $\{\frac{\partial }{\partial z_{i}}\}(dz_{i})$ is a local framing for $%
T^{(1,0)}\left( X)\right) (\bigwedge^{(1,0)}(X))$. Formulas \ref{higgsform}
and \ref{higgsform2} imply that $\theta $ is actually a holomorphic section
of $Hom(E)\otimes \Lambda ^{1,0}(X)$ and the condition $\left[ \theta
(Z),\theta (W)\right] =0$ follows from the assumption that k is odd as
follows:

\begin{equation*}
\lbrack \theta (Z),\theta (W)](\phi )=\{i(Z)\varpi \}\wedge \{i(W)\varpi
\}\wedge \phi -\{i(W)\varpi \}\wedge \{i(Z)\varpi \}\wedge \phi =0
\end{equation*}
because $i(Z)\varpi $ is a form of even degree. This same idea shows that if 
$\varpi $ is a sum of holomorphic forms of possibly different odd degrees,
then \ref{higgsform} also defines a Higgs structure on $E$. If $\varpi $ is
a holomorphic k-form, where $k$ is not necessarily assumed to be odd, then a
``super Higgs'' structure can be defined on $E$ if we define a new bracket
operation $``\ [\ ,\ ]^{\varpi }\ "$ in $Hom(E)$ by the prescription 
\begin{equation*}
\lbrack A\ ,\ B]^{\varpi }\ (\phi )=(AB-(-1^{\deg (\varpi )})BA)(\phi )
\end{equation*}

We now examine some examples of this Higgs form for specific values of $k$
in a purely linear algebraic setting . Let ($V_{\mathbb{R}},J)$ be a real
vector space with a complex structure $J$, $V=V_{\mathbb{R}}\bigotimes 
\mathbb{C=}V^{(1,0)}\bigoplus V^{(0,1)}$ be the complexification and
decomposition into $\pm i$ eigenspaces of $J$. Let $\varpi \in \Lambda
^{(k,0)}(V)$ ($k$ odd) and define $\theta \in Hom(\bigoplus_{p\geq 0}\Lambda
^{(p,0)}(V))\bigotimes \Lambda ^{(1,0)}(V)$ by \ref{higgsform}. If $k=1,$
then \ref{higgsform} yields $\theta (\phi )=\phi \bigotimes \varpi $, that
is, $\theta (Z)(\phi )=\varpi (Z)\phi .$ If $k=n$ is odd, then 
\begin{equation*}
\theta (\phi )=\left\{ 
\begin{array}{c}
\phi i(\bullet )\varpi \text{ if }\deg (\phi )=0 \\ 
\varpi \bigotimes \phi \text{ if }\deg (\phi )=1 \\ 
0\text{ if }\deg (\phi )\geq 2
\end{array}
\right.
\end{equation*}

the middle expression means $\theta (Z)(\phi )=\phi (Z)\varpi $ and these
formulas follow from $(i(Z)\varpi )\wedge \phi =i(Z)(\varpi \wedge \phi
)+\varpi \wedge (i(Z)\phi )$, which is valid for any form $\phi $. These
examples show the kernel of $\theta $ is $0$ if $\varpi \neq 0$ in the
interesting cases where $\theta $ could act nontrivially. In general we have

\begin{proposition}
\label{kertheta} i. For $\phi \in \bigoplus_{p\geq 0}\Lambda ^{(p,0)}(V)$, $%
\theta (\phi )=0$ $\Longleftrightarrow $ $\varpi \wedge \phi =0$ and $\varpi
\wedge i(Z)\phi =0$ $\forall $ $Z\in V.$\newline
\ ii. Let $h$ be any Hermitian metric on $\bigoplus_{p\geq 0}\Lambda
^{(p,0)}(V)$, and let $\overline{\theta _{h}}$ be the $h$-adjoint of $\theta
,$ $h(\overline{\theta _{h}}(Y)\ \phi ,\psi )=h(\phi ,\theta (\ \overline{Y}%
\ )\ \psi )$. Then $\overline{\theta _{h}}(\overline{Z})\ \phi =0\;\forall Z$
$\in V\Longleftrightarrow \left( \varepsilon (\varpi \right) )^{\ast
_{h}}\phi =0$ and $\left( \varepsilon (\varpi \right) )^{\ast
_{h}}i(Z)^{\ast _{h}}\phi =0\;\forall Z\in V$ where $^{\ast _{h}}$ means
adjoint with respect to $h.$
\end{proposition}

\begin{proof}
i. We need only consider $Z\in V^{(1,0)}$. The formula $(i(Z)\varpi )\wedge
\phi =i(Z)(\varpi \wedge \phi )+\varpi \wedge i(Z)\phi $ makes $%
\Longleftarrow $ clear. If $(i(Z)\varpi )\wedge \phi =0$ $\forall Z$ then $%
((\varepsilon (\eta )i(Z)\varpi )\wedge \phi =0$ $\forall \eta \in \Lambda
^{(1,0)}(V)$.Therefore $0=\sum_{j}(\varepsilon (Z_{j}^{\ast })i(Z_{j})\varpi
)\wedge \phi $ where $\{Z_{j}\}(\{Z_{j}^{\ast }\})$ is a basis (dual) for $%
V^{(1,0)}(\Lambda ^{(1,0)}(V))$ but this sum also equals $k$ $\varpi \wedge
\phi $ due to the identity $\sum_{j}\varepsilon (Z_{j}^{\ast
})i(Z_{j})\varpi =k$ $\varpi $ which is valid for any $(k,0)$ form
(seemingly most easily proved by computing on basis elements $%
Z_{j_{1}}^{\ast }\wedge Z_{j_{2}}^{\ast }\wedge \cdots \wedge
Z_{j_{k}}^{\ast }$). Thus $\theta (\phi )=0$ implies $\varpi \wedge \phi =0.$
Therefore the assumption in $\Longrightarrow $ yields $0=(i(Z)\varpi )\wedge
\phi =i(Z)(\varpi \wedge \phi )+\varpi \wedge i(Z)\phi =\varpi \wedge
i(Z)\phi $ $\forall Z$. \newline
ii. $\overline{\theta _{h}}(\overline{Z})\ \phi =0\;\forall Z$ $\in
V\Longleftrightarrow h(\phi ,(i(\ Z\ )\varpi )\wedge \ \psi )=0\;\forall Z$
and $\forall \psi $. Replacing $\psi $ with $\varepsilon (Z^{\ast })\psi $,
this implies $h(\phi ,\sum_{j}\varepsilon (Z_{j}^{\ast })(i(Z_{j})\varpi
)\wedge \psi )=0,$ so $h(\phi ,\varpi \wedge \psi )=0$ i.e. $h(\varepsilon
(\varpi )^{\ast _{h}}\phi ,\psi )=0,$and thus $\varepsilon (\varpi )^{\ast
_{h}}\phi =0.$ The rest is as is in part i.
\end{proof}

Let us call a (positive definite) Hermitian metric $h$ on $%
\bigoplus_{p,q\geq 0}\Lambda ^{(p,q)}(V)$ \textit{standard }if $h$ is the
unique extension to $\bigoplus_{p,q\geq 0}\Lambda ^{(p,q)}(V)$ of a (real)
metric on $V_{\mathbb{R}}$ for which $J$ is orthogonal such that $\Lambda
^{(p,q)}(V)$ is orthogonal to $\Lambda ^{(p^{\prime },q^{\prime })}(V)$ if $%
(p,q)\neq (p^{\prime }q^{\prime })$ and $Z_{i_{1}}^{\ast }\wedge
Z_{i_{2}}^{\ast }\wedge \cdots \wedge Z_{i_{p}}^{\ast }\wedge \overline{%
Z_{j_{1}}^{\ast }}\wedge \cdots \wedge \overline{Z_{j_{q}}^{\ast }}$, $1\leq
i_{1}<\cdots <i_{p}\leq n$, $1\leq j_{1}<\cdots <j_{q}\leq n$ is an
orthonormal basis for $\Lambda ^{(p,q)}(V)$ if $\{Z_{j}\}_{j=1}^{n}$ is an
orthonormal basis for $V^{(1,0)}.$ If $h$ is standard\ then one has the
usual isomorphisms $\#:V^{\ast }\rightarrow V$ and $\flat :V\rightarrow
V^{\ast }$ and then $\varepsilon (\varpi )^{\ast h}=i(\varpi ^{\#})$ and $%
i(Z)^{\ast h}=\varepsilon (Z^{\flat })$. One proves the following statement

\begin{equation}
\text{For a standard }h\text{, }\overline{\theta _{h}}(\phi
)=0\Longleftrightarrow i(\varpi ^{\#})\phi =0\text{ and }\varpi \wedge
\varepsilon (Z^{\flat })\phi =0\text{\ }\forall Z\in V
\end{equation}
\ Let us call a (positive definite) Hermitian metric $h$ on $%
\bigoplus_{p\geq 0}\Lambda ^{(p,0)}(V)$ \textit{natural }if $h$ makes $%
\Lambda ^{(p,0)}(V)$ orthogonal to $\Lambda ^{(p^{\prime },0)}(V)$ if $\
p\neq p^{\prime }.$ One verifies that 
\begin{equation}
h\text{ natural on}\bigoplus_{p\geq 0}\Lambda ^{(p,0)}(V)\Longrightarrow 
\overline{\theta _{h}}:\Lambda ^{(a,0)}(V)\rightarrow \Lambda
^{(a-k+1,0)}(V)\bigotimes \Lambda ^{(1,0)}(V)  \label{hnatural}
\end{equation}

For use later in giving examples of K\"{a}hler manifolds which admit
Higgs-Hermitian-Yang-Mills metrics (\ref{HHYMex}) we now give a formula for
the linear algebraic operator $T_{h}(s)$ defined by 
\begin{eqnarray}
T_{h}(s) &:&=\overset{n}{\underset{i=1}{\sum }}[\theta (Z_{i})\overline{%
,\theta _{h}}(\overline{Z_{i}})]s  \notag \\
h(T_{h}(s),s) &=&\overset{n}{\underset{i=1}{\sum }\{}\Vert \overline{\theta
_{h}}(\overline{Z}_{i})s\Vert _{h}^{2}-\Vert \theta (Z_{i})s\Vert _{h}^{2}\}
\label{T-op}
\end{eqnarray}
\ (cf. \ref{boch5}), in the case where the Hermitian metric $h$ on $%
\bigoplus_{p,q\geq 0}\Lambda ^{(p,q)}(V)$ is standard, $s\in
\bigoplus_{p\geq 0}\Lambda ^{(p,0)}(V)$ and the $\theta $ operator is
defined using an element $\varpi =aZ_{_{1}}^{\ast }\wedge Z_{_{2}}^{\ast
}\wedge \cdots \wedge Z_{n}^{\ast }\in \Lambda ^{(n,0)}(V)$ ($n$ odd, 
\TEXTsymbol{>} 1) where $\{Z_{j}\}_{j=1}^{n}$ is an orthonormal basis for $%
V^{(1,0)}$. In this setting, one has the identity $i((Z_{i_{1}}^{\ast
}\wedge Z_{i_{2}}^{\ast }\wedge \cdots \wedge Z_{i_{p}}^{\ast
})^{\#})=i(Z_{i_{p}})i(Z_{i_{p-1}})\cdots i(Z_{i_{1}})$ for any $p$. If $%
\deg s\geq 2$ then $\theta (Z_{i})s=\{i(Z_{i})\varpi \}\wedge s=0$ and if $%
\deg s\leq n-2$ then $\overline{\theta _{h}}(\overline{Z}_{i})s=i(\pm
(aZ_{_{1}}^{\ast }\wedge Z_{_{2}}^{\ast }\wedge \cdots \wedge \widehat{%
Z_{_{i}}^{\ast }}\wedge \cdots \wedge Z_{n}^{\ast })^{\#})s=0$. Thus $%
T_{h}(s)=0$ if $2\leq \deg s\leq n-2.$ It is straightforward to check that $%
h(T_{h}(s),s)=f[i]\left\| \varpi \right\| _{h}^{2}\left\| s\right\| _{h}^{2}$%
, where $s\in \Lambda ^{(i,0)}(V)$, with $f[i]=0$, if $2\leq i\leq n-2$, $%
f[0]=-n,f[1]=-1,f[n]=n$ and (because $T_{h}$ must have trace 0, or by
similar computations) $f[n-1]=1.$ By polarizing, we get 
\begin{gather}
\varpi \in \Lambda ^{(n,0)}(V),s\in \Lambda ^{(i,0)}(V)\Rightarrow  \notag \\
T_{h}(s)=\left\| \varpi \right\| _{h}^{2}f[i]s  \label{f[i] values} \\
f[i]=\left\{ 
\begin{array}{c}
-n\text{ if }i=0 \\ 
-1\text{ if }i=1 \\ 
0\text{ if }2\leq i\leq n-2 \\ 
1\text{ if }i=n-1 \\ 
n\text{ if }i=n
\end{array}
\right.  \notag
\end{gather}

We remark that one can prove the following identity: if $(V,h)$ are as
above, but now $\varpi =aZ_{i_{1}}^{\ast }\wedge Z_{i_{2}}^{\ast }\wedge
\cdots \wedge Z_{i_{k}}^{\ast }\in \Lambda ^{(k,0)}(V)$ is a simple $(k,0)-$%
form ($k$ odd), then 
\begin{gather*}
T_{h}(s)= \\
-\left\| \varpi \right\| _{h}^{2}\{ks+ \\
\sum_{r=1}^{\min (k-2,\deg s)}(k-r)(-1)^{r}\sum_{1\leq t_{1}<t_{2}<\cdots
t_{r}\leq k}\varepsilon (Z_{i_{t_{1}}}^{\ast }\wedge Z_{i_{t_{2}}}^{\ast
}\wedge \cdots \wedge Z_{i_{t_{r}}}^{\ast })i((Z_{i_{t_{1}}}^{\ast }\wedge
Z_{_{i_{t_{2}}}}^{\ast }\wedge \cdots \wedge Z_{i_{t_{r}}}^{\ast
})^{^{\#}})s\}
\end{gather*}

Consequently if $s\in \Lambda ^{(i,0)}((span\{Z_{i_{1}}^{\ast
},Z_{i_{2}}^{\ast },\cdots ,Z_{i_{k}}^{\ast }\}^{\bot })$ , then $%
T_{h}(s)=-k\left\| \varpi \right\| _{h}^{2}s$, while if $s\in \Lambda
^{(i,0)}((span\{Z_{i_{1}}^{\ast },Z_{i_{2}}^{\ast },\cdots ,Z_{i_{k}}^{\ast
}\})$, then one can show $T_{h}(s)=\left\| \varpi \right\| _{h}^{2}F[i]s$,
with $F[0]=-k,F[1]=-1,F[k-1]=1,F[0]=k,F[i]=0$ if $2\leq i\leq k-2$.

Now consider again the differential geometric setting described in the
beginning of section 2. $E\rightarrow X$ is the holomorphic vector bundle $%
E=\bigoplus_{p\geq 0}\Lambda ^{(p,0)}(X)$, $\varpi \in \vartheta \mathrm{%
\Gamma }\Lambda ^{(k,0)}(X)$,and $\theta $ the Higgs form defined by \ref
{higgsform}. Let $h$ be any Hermitian metric on $E$ and let $g$ be any
Hermitian metric on $TX\bigotimes \mathbb{C}$ (we do not assume any a priori
relation between $g$ and $h$). In this case formula \ref{bochner} becomes 
\begin{equation}
h(F_{h}(Z,\overline{Z})s,s)=h(\Theta (Z,\overline{Z})s,s)+\Vert (\varepsilon
(i(Z)\varpi ))^{\ast _{h}}s\Vert _{h}^{2}-\Vert (i(Z)\varpi )\wedge s\Vert
_{h}^{2}  \label{boch3}
\end{equation}

\begin{remark}
1.If $k=\deg \varpi =1$, then $(i(Z)\varpi )\wedge s=(i(Z)\varpi
)s,(\varepsilon (i(Z)\varpi )^{\ast _{h}}s=\overline{(i(Z)\varpi )}s$ (even
if $h$ is not natural) and then \ref{boch3} becomes $h(F_{h}(Z,\overline{Z}%
)s,s)=h(\Theta (Z,\overline{Z})s,s)$. In fact, the operator corresponding to 
$\theta \wedge \overline{\theta _{h}}+\overline{\theta _{h}}\wedge \theta $
(cf. \ref{F11}) is zero. Partly this reason we will assume $k\geq 3$ unless
specified otherwise. Another reason for assuming $k\geq 3$ is that we want
to consider solutions to the equation $(i(Z)\varpi )\wedge s=0$ $\forall Z$
(locally defined) holomorphic tangent vector fields, and $s\in \vartheta
\Gamma E$. If $\varpi $ is a $1-form$, then this would imply either $\varpi
=0$ or $s=0$.\newline
2. Note that $F_{h}$ does not annihilate functions on\ $X,$ i.e. sections of 
$\Lambda ^{(0,0)}(X),$ unlike $\Theta $. In particular, we conclude for the
constant section $1\in \Gamma \Lambda ^{(0,0)}(X)$, that $h(F_{h}(Z,%
\overline{Z})1,1)=\Vert (\varepsilon (i(Z)\varpi ))^{\ast _{h}}1\Vert
_{h}^{2}-\Vert (i(Z)\varpi )\Vert _{h}^{2}$ and if $h$ is natural, then $%
h(F_{h}(Z,\overline{Z})1,1)=-\Vert (i(Z)\varpi )\Vert _{h}^{2}$
\end{remark}

If the $\deg s\geq n-k+2$, $(i(Z)\varpi )\wedge s=0$. If $h$ is a natural
metric then $\deg (\varepsilon (i(Z)\varpi )^{\ast _{h}}s$ $=\deg s-k+1$ and
hence if $\deg s\leq k-2$ then $(\varepsilon (i(Z)\varpi ))^{\ast _{h}}s=0.$
Now $k-2\geq n-k+2$ $\Longleftrightarrow k\geq \frac{n}{2}+2$, so for this
range of $k$ both of the last two terms on the right hand side of \ref{boch3}
vanish. We summarize these observations below. 
\begin{gather}
h(F_{h}(Z,\overline{Z})s,s)\geq h(\Theta (Z,\overline{Z})s,s)\text{ }\forall
s\in C^{\infty }\Gamma \bigoplus_{a\geq n-k+2}\Lambda ^{(a,0)}(X)
\label{ineq1} \\
h\,natural\Rightarrow h(F_{h}(Z,\overline{Z})s,s)\leq h(\Theta (Z,\overline{Z%
})s,s)\text{ }\forall s\in C^{\infty }\Gamma \bigoplus_{a\leq k-2}\Lambda
^{(a,0)}(X)\text{\label{ineq2}} \\
h\,natural,k\geq \frac{n}{2}+2\text{ }\Rightarrow \text{ }F_{h}(Z,\overline{Z%
})s=\Theta (Z,\overline{Z})s\text{ }\forall s\in C^{\infty }\Gamma
\bigoplus_{n-k+2\leq a\leq k-2}\Lambda ^{(a,0)}(X)  \notag
\end{gather}

\begin{theorem}
\label{thmw=0}Assume $(X,g)$ is a compact Hermitian manifold of complex
dimension n, $E$ $\rightarrow X$ is the Higgs bundle given by \ref{higgsform}
and $h$ is a Hermitian metric on $E$. If for all sections $t$ of $E$, $0\geq
h(i\Lambda F_{h}t,t)$, pointwise$,$ then 
\begin{equation*}
s\in \vartheta \Gamma \bigoplus_{a\geq n-k+2}\Lambda
^{(a,0)}(X)\Longrightarrow \nabla ^{h}s\equiv 0,(\varepsilon (i(Z)\varpi
)^{\ast _{h}}s\equiv 0\forall Z,and\text{ }i\Lambda F_{h}(s)\equiv 0
\end{equation*}
If $0\geq h(i\Lambda F_{h}t,t)$ for all sections $t$ of $E$ and $k\geq \frac{%
n}{2}+1,$ then $\varpi =0$.
\end{theorem}

\begin{proof}
\ref{boch2} in this setting can be written 
\begin{equation}
\overset{}{h(i\Lambda F_{h}s,s)-h(i\Lambda \Theta s,s)\underset{}{=}\sum \{}%
\Vert (\varepsilon (i(Z_{i})\varpi ))^{\ast _{h}}s\Vert _{h}^{2}-\Vert
\{i(Z_{i})\varpi \}\wedge s\Vert _{h}^{2}\}  \label{boch5}
\end{equation}
The argument of Lemma \ref{lemma1} implies the result in the first line of
the theorem, because $\deg (s)\geq n-k+2$ implies $\{i(Z)\varpi \}\wedge s=0$
$\forall Z$. To prove the second statement, note that $k\geq \frac{n}{2}+1$
implies that $\deg (\varpi )=k\geq n-k+2$ so we can use the first argument
to conclude that ($\varpi $ is $h$-parallel and ) $(\varepsilon (i(Z)\varpi
))^{\ast _{h}}\varpi =0\forall Z.$ From the second part of Proposition \ref
{kertheta} we get $(\varepsilon (\varpi ))^{\ast _{h}}\varpi =0$. This
yields $h((\varepsilon (\varpi ))^{\ast _{h}}\varpi ,t)=0$ for all sections $%
t$ of $E$, and taking $t=1$ implies $\Vert \varpi \Vert _{h}^{2}=0$.
\end{proof}

If $i\Lambda F_{h}$ is quasinegative and $s$ is as in Theorem \ref{thmw=0}
then $s=0$.

In the context of general Higgs bundles, the vanishing of the $F_{h}^{(2,0)}$
and $F_{h}^{(0,2)}$ is equivalent to the Higgs form being parallel (cf.
Proposition \ref{prop}). The next result examines the case of those Higgs
bundles defined by \ref{higgsform}, and with a special metric. The result
will be used at the end of section 3 in a vanishing theorem for $(p,q)$
forms with values in $E$.

\begin{proposition}
\label{wparallel}Let $(X,g)$ \ be a K\"{a}hler manifold, and extend $g$ to a
standard metric on $\Lambda ^{\ast }(X)\bigotimes \mathbb{C}$. In the above
notation let the metric $h$ on $E:=\overset{n}{\underset{p=0}{\bigoplus }}%
\bigwedge^{(p,0)}(X)$ be $h=g$ $\ $Then $F_{h}^{(2,0)}=0$ $%
\Longleftrightarrow \nabla \varpi =0$ ($\Leftrightarrow F_{h}^{(0,2)}=0$ by
Proposition\ref{prop}).
\end{proposition}

\begin{proof}
From Proposition \ref{prop} it follows that on any Higgs bundle $(E,\theta )$
with a Hermitian metric $h$, $F_{h}^{(2,0)}=0$ at a point $p\Leftrightarrow $
$\partial \theta _{s}^{t}=0$ at $p$ where $\theta e_{\alpha }=\sum_{\beta
}e_{\beta }\bigotimes \theta _{\alpha }^{\beta }$ for a local frame $%
\{e_{\alpha }\}_{\alpha =1}^{r}$ of $E$ adapted to $p$. In the case we are
considering, let $\{e_{\alpha }\}_{\alpha =1}^{r}$ ($r=2^{n}$) be a local
frame of $E=\bigoplus_{p=0}^{n}\Lambda ^{(p,0)}(X)$ $h$-adapted to $p$. Then
throughout the neighborhood where $\{e_{\alpha }\}_{\alpha =1}^{r}$ is
defined we can write $\theta e_{\alpha }=\sum_{i=1}^{n}(\{i(\frac{\partial }{%
\partial z_{i}})\varpi \}\wedge e_{\alpha })\bigotimes dz_{i}=\sum_{\beta
=1}^{r}(e_{\beta }\bigotimes \sum_{i=1}^{n}A_{i,\alpha }^{\beta }dz_{i})$
where $\{i(\frac{\partial }{\partial z_{i}})\varpi \}\wedge e_{\alpha
}=\sum_{\beta =1}^{r}A_{i,\alpha }^{\beta }e_{\beta }$ . Therefore we have $%
\theta _{\alpha }^{\beta }=\sum_{i=1}^{n}A_{i,\alpha }^{\beta }dz_{i}$ and $%
\partial \theta _{\alpha }^{\beta }=\sum_{i=1}^{n}\partial A_{i,\alpha
}^{\beta }\wedge dz_{i}$. One computes that $A_{i,\alpha }^{\beta
}=\sum_{\gamma =1}^{r}h^{\beta \overline{\gamma }}h(\{i(\frac{\partial }{%
\partial z_{i}})\varpi \}\wedge e_{\alpha },e_{\gamma })$,where $h_{a%
\overline{b}}=h(e_{a},e_{b})$ and $(h^{s\overline{t}})=(h_{a\overline{b}%
})^{-1}$. Up to this point we have not used the assumption that $h$ is a
K\"{a}hler metric. We now exploit this assumption by writing $\partial
=\sum_{j}\varepsilon (\frac{\partial }{\partial z_{j}})\nabla _{\frac{%
\partial }{\partial z_{j}}}$ where $\{\frac{\partial }{\partial z_{j}}%
\}_{j=1}^{n}$ is a local holomorphic frame for $T^{(1,0)}(X)$ which is also $%
h-$adapted to $p$. Then, at $p$, the following equalities hold: 
\begin{eqnarray*}
\partial \theta _{\alpha }^{\beta } &=&\sum_{i,j}\frac{\partial }{\partial
z_{j}}h(\{i(\frac{\partial }{\partial z_{i}})\varpi \}\wedge e_{\alpha
},e_{\beta })dz_{j}\wedge dz_{i} \\
&=&\sum_{i,j}h(\nabla _{\frac{\partial }{\partial z_{j}}}(\{i(\frac{\partial 
}{\partial z_{i}})\varpi \}\wedge e_{\alpha }),e_{\beta })dz_{j}\wedge dz_{i}
\\
&=&\sum_{i,j}h((\nabla _{\frac{\partial }{\partial z_{j}}}\{i(\frac{\partial 
}{\partial z_{i}})\varpi \})\wedge e_{\alpha },e_{\beta })dz_{j}\wedge dz_{i}
\\
&=&\sum_{i,j}h((i(\frac{\partial }{\partial z_{i}})\nabla _{\frac{\partial }{%
\partial z_{j}}}\varpi )\wedge e_{\alpha },e_{\beta })dz_{j}\wedge dz_{i} \\
&=&\sum_{i<j}h([i(\frac{\partial }{\partial z_{i}})\nabla _{\frac{\partial }{%
\partial z_{j}}}\varpi -i(\frac{\partial }{\partial z_{j}})\nabla _{\frac{%
\partial }{\partial z_{i}}}\varpi ]\wedge e_{\alpha },e_{\beta
})dz_{j}\wedge dz_{i}
\end{eqnarray*}
We conclude: $\partial \theta _{\alpha }^{\beta }(p)=0\Leftrightarrow h([i(%
\frac{\partial }{\partial z_{i}})\nabla _{\frac{\partial }{\partial z_{j}}%
}\varpi -i(\frac{\partial }{\partial z_{j}})\nabla _{\frac{\partial }{%
\partial z_{i}}}\varpi ]\wedge e_{\alpha },e_{\beta })(p)=0$ $\forall i<j$ .
Thus $\nabla \varpi =0\Rightarrow $ $\partial \theta _{\alpha }^{\beta }=0$
and hence that $F_{h}^{(2,0)}=0.$ Conversely, $F_{h}^{(2,0)}=0\Rightarrow
\partial \theta _{\alpha }^{\beta }(p)=0\Rightarrow h([i(\frac{\partial }{%
\partial z_{i}})\nabla _{\frac{\partial }{\partial z_{j}}}\varpi -i(\frac{%
\partial }{\partial z_{j}})\nabla _{\frac{\partial }{\partial z_{i}}}\varpi
]\wedge e_{\alpha },e_{\beta })(p)=0$ for any adapted frame $\{e_{\alpha
}\}_{\alpha =1}^{r}$ and we can conclude that $F_{h}^{(2,0)}=0\Rightarrow (i(%
\frac{\partial }{\partial z_{i}})\nabla _{\frac{\partial }{\partial z_{j}}%
}\varpi -i(\frac{\partial }{\partial z_{j}})\nabla _{\frac{\partial }{%
\partial z_{i}}}\varpi )(p)=0$ $\forall i,j$ e.g. by choosing $e_{\alpha
}=1\in \Lambda ^{(0,0)}(X)\subset E.$ Finally we have 
\begin{gather*}
i(\frac{\partial }{\partial z_{i}})\nabla _{\frac{\partial }{\partial z_{j}}%
}\varpi =i(\frac{\partial }{\partial z_{j}})\nabla _{\frac{\partial }{%
\partial z_{i}}}\varpi \Rightarrow \\
\sum_{i}\varepsilon (dz_{i})i(\frac{\partial }{\partial z_{i}})\nabla _{%
\frac{\partial }{\partial z_{j}}}\varpi =\sum_{i}\varepsilon (dz_{i})i(\frac{%
\partial }{\partial z_{j}})\nabla _{\frac{\partial }{\partial z_{i}}}\varpi
\\
\Rightarrow k\text{ }\nabla _{\frac{\partial }{\partial z_{j}}}\varpi
=\sum_{i}\varepsilon (dz_{i})i(\frac{\partial }{\partial z_{j}})\nabla _{%
\frac{\partial }{\partial z_{i}}}\varpi \Rightarrow \\
k\text{ }\sum_{j}\varepsilon (dz_{j})\nabla _{\frac{\partial }{\partial z_{j}%
}}\varpi =\sum_{i,j}\varepsilon (dz_{j})\varepsilon (dz_{i})i(\frac{\partial 
}{\partial z_{j}})\nabla _{\frac{\partial }{\partial z_{i}}}\varpi
\Rightarrow \\
k\text{ }\partial \varpi =-\sum_{i,j}\varepsilon (dz_{i})\varepsilon
(dz_{j})i(\frac{\partial }{\partial z_{j}})\nabla _{\frac{\partial }{%
\partial z_{i}}}\varpi \\
\Rightarrow k\text{ }\partial \varpi =-k\text{ }\partial \varpi \Rightarrow
\partial \varpi =0.
\end{gather*}
From the third line above we also get $k\nabla _{\frac{\partial }{\partial
z_{j}}}\varpi =\nabla _{\frac{\partial }{\partial z_{j}}}\varpi -\sum_{i}i(%
\frac{\partial }{\partial z_{j}})\varepsilon (dz_{i})\nabla _{\frac{\partial 
}{\partial z_{i}}}\varpi =\nabla _{\frac{\partial }{\partial z_{j}}}\varpi $ 
$-i(\frac{\partial }{\partial z_{j}})\partial \varpi =\nabla _{\frac{%
\partial }{\partial z_{j}}}\varpi $ because $\partial \varpi =0$ . Finally
we conclude $(k-1)\nabla _{\frac{\partial }{\partial z_{j}}}\varpi =0,$and
hence $\nabla \varpi =0.$
\end{proof}

In preparation for the examination of stability questions for the Higgs
bundle $E,$ we consider general Higgs subbundles of $E.$ Suppose $P\subset E$
is a Higgs subbundle of $E.$ This means that if $s$ is a local section of $P$%
, then $(i(Z)\varpi )\wedge s$ is also a local section of $P$. If $h$ is a
Hermitian metric on $P$, then \ref{boch3} and \ref{boch5} apply to the Higgs
bundle $P$ with the Hermitian metric $h$. Additionally, the proof of Theorem 
\ref{thmw=0} works as well in this setting, which we include as a :

\begin{remark}
\label{rmksbbdl}Let $P\subset E$ be a Higgs subbundle and let $h$ be a
Hermitian metric on $P$, so all $h$ Hermitian data applies to $P$. If $0\geq
i\Lambda F_{h}$, i.e. $i\Lambda F_{h}$ is a pointwise negative semidefinite
operator, then any holomorphic section $s$ of $P$ which is a $(p,0)$-form
with $p\geq n-k+2$, (or a sum of such forms) must be parallel for the
Hermitian connection of $h$.\ If $i\Lambda F_{h}$ is quasinegative, then any
such $s$ must be $0$.
\end{remark}

We now examine the question of stability for the Higgs bundles defined by 
\ref{higgsform}. Assume ($X,g)$ is a compact K\"{a}hler manifold. If $%
E\rightarrow X$ is any holomorphic vector bundle over $X$, then $E$ is said
to be \textit{stable (semistable)} (\cite{kobayashi},Chapter V, sections
5-7) if for every nontrivial coherent analytic subsheaf $\frak{F}$ of the
sheaf $\vartheta (E)$ of germs of holomorphic sections of $E$ the following
inequality holds: 
\begin{equation}
\mu (\frak{F})<(\leq )\mu (E)  \label{stable}
\end{equation}
If $F\subset E$ is any holomorphic subbundle of $E,$ $\mu (F)$ (the
``slope'' of $F$ ) is defined to be $\mu (F)=\frac{\deg (F)}{rank(F)}=\frac{%
\int_{X}c_{1}(F)\wedge \omega ^{n-1}}{rank(F)},\omega $ being the K\"{a}hler
form of $g$ and $c_{1}(F)$ the first Chern class of $F.$ If $\frak{F}$ is
the sheaf $\vartheta (F)$ of germs of holomorphic sections of $F,\mu (\frak{F%
})$ means $\mu (F).$ A coherent subsheaf of $\vartheta (E)$ need not arise
as the sheaf of germs of such a subbundle $F$, i.e. $\frak{F}$ need not be
locally free. Nevertheless, there is a well-defined rank for $\frak{F}$,
because $\frak{F}$ is locally free outside a set of codimension at least 2.
There is also a holomorphic line bundle associated to $\frak{F}$ ,``$\det (%
\frak{F})$'' and one defines $\ c_{1}(\frak{F})$ to be $c_{1}(\det (\frak{F}%
))$ and $\mu (\frak{F})=\frac{\int_{X}c_{1}(\frak{F})\wedge \omega ^{n-1}}{%
rank(\frak{F})}$. In case $\frak{F}$ arises from a vector bundle $F$ these
definitions agree with the standard vector bundle ones.

A Hermitian metric $h$ on the holomorphic vector bundle $E$ over $(X,g)$ is
said to be an Einstein-Hermitian metric (\cite{kobayashi}, chapter IV) or a
Hermitian-Yang-Mills (HYM) metric (\cite{yau}) if $i\Lambda \Theta =cId_{E}$%
, $\Theta $ being the Hermitian curvature of $h$, and where $c$ is a
constant determined by the rank and degree of $E$ and the (class of ) the
K\"{a}hler form of $g$ (cf. \cite{kobayashi} chapter IV, section 2). If $E$
admits such a metric $h$, then $E$ is semistable and splits into a direct
sum of holomorphic, stable subbundles with the same slope (\cite{kobayashi}
chapter V section 8, Theorem 8.3). The converse theorem conjectured by
Kobayashi was proved in \cite{yau}: A stable holomorphic vector bundle over
a compact K\"{a}hler manifold admits a unique Hermitian-Yang-Mills metric.

In the category of Higgs bundles over compact K\"{a}hler manifolds, $%
(E,\theta )\rightarrow (X,g)$, $E$ is said to be \textit{Higgs} \textit{%
stable (Higgs semistable)} (\cite{simpConstVHS}) if for every nontrivial
coherent analytic subsheaf $\frak{F}$ satisfying $\theta :\frak{F\rightarrow
F}\bigotimes \vartheta (\Lambda ^{(1,0)}(X))$ (i.e. a Higgs subsheaf) the
inequality in \ref{stable} holds. A Hermitian metric $h$ on the Higgs bundle 
$(E,\theta )\rightarrow (X,g)$ is said to be a (Higgs-)Hermitian-Yang-Mills
(HHYM) metric (\cite{simpIHES}) if $i\Lambda F_{h}=cId_{E}$, where $F_{h}$
is defined in \ref{F}. . Again it is true that if $(E,\theta )$ admits such
a metric $h$, then $(E,\theta )$ is (Higgs)semistable and splits into a
direct sum of holomorphic, (Higgs)stable subbundles with the same slope
(``polystability''\cite{simpIHES},Theorem 1, page 19), because the proof of 
\cite{kobayashi} chapter V section 8, Theorem 8.3 can be modified for the
Higgs category, and the inequalities still go the right way. The converse of
this theorem, for compact and certain classes of noncompact K\"{a}hler
manifolds, is due to Simpson (\cite{simpConstVHS}, see also \cite{simpIHES})
and plays an important part in the results described at the beginning of
section 1.

One would like to know when an HHYM $h$ exists for the Higgs bundles defined
by (\ref{higgsform}) for a $X$ a compact K\"{a}hler manifold. The results we
present below (\ref{bigthm}) indicate that such metrics may be quite rare
for such $X$. In order to get some information about such metrics we give
examples of HHYM metrics in noncompact cases where there are no topological
or complex-analytic obstructions to their existence.

Let $(X,g)$ be complex n-dimensional with K\"{a}hler metric $g.$ Assume the
following properties are satisfied:

i. \ $g$ is K\"{a}hler-Einstein

ii. given any constant $C$, there is a smooth function $f:X\rightarrow 
\mathbb{C}$ such that $\overline{\square }_{g}(f)=C$

iii. $\varpi \in \Lambda ^{(n,0)}(X)$ is a holomorphic $n-$form with
constant $g-$length

Then the Higgs bundle $(E=\bigoplus_{p\geq 0}\Lambda ^{(p,0)}(X),\theta )$
admits\ a HHYM metric $g^{\prime }$, $i\Lambda F_{g^{\prime }}=cId_{E}$ with
any number $c$ (note that the $\Lambda $ in $i\Lambda F_{g^{\prime }}$
refers to interior multiplication by the $g-$dual to the $g-$K\"{a}hler
form, we use $g$ for all Riemannian data on $X$). In fact we will now show
that such a $g^{\prime }$\ can be obtained by taking the standard extension
of $g$ to $E$ and changing it conformally on each $\Lambda ^{(p,0)}(X)$
(with a conformal factor depending on $p$). $X=\mathbb{C}^{n}$ with the
standard metric and a constant coefficient $(n,0)-$form $\varpi $ is of
course an example of such a manifold, and $f(z_{1},z_{2},...,z_{n})=\frac{C}{%
n}\sum_{i}\left| z_{i}\right| ^{2}+$any harmonic function yields condition
ii. That condition excludes the possibility of $X$ being compact.

For $(X,g)$ satisfying i, ii and iii, say that $Ric(g,T^{(1,0)}(X))=-\lambda
Id_{T^{(1,0)}(X)},$ so for the induced action of the Ricci curvature on $%
\Lambda ^{(1,0)}(X),$ $Ric(g,\Lambda ^{(1,0)}(X))=\lambda Id_{\Lambda
^{(1,0)}(X)}$. Extend $g$ as a standard Hermitian metric to $E.$ Then the
curvature $\Theta _{g}$ of the corresponding Hermitian connection on $E$
then satisfies $i\Lambda \Theta _{g}=\binom{n-1}{p-1}\lambda Id_{\Lambda
^{(p,0)}(X)}$ on $\Lambda ^{(p,0)}(X)$. The decomposition $%
E=\bigoplus_{p\geq 0}\Lambda ^{(p,0)}(X)$ is $g-$orthogonal, hence each
factor $\Lambda ^{(p,0)}(X)$ is invariant under the Hermitian connection $%
\nabla ^{g},$ i.e. each factor is totally geodesic. Thus $i\Lambda \Theta
_{(\Lambda ^{(p,0)}(X),g\mid \Lambda ^{(p,0)}(X))}=i\Lambda \Theta _{g}\mid
_{\Lambda ^{(p,0)}(X)}=\binom{n-1}{p-1}\lambda Id_{\Lambda ^{(p,0)}(X)}$.
Let $g^{\prime }$ be the Hermitian metric on $E$ uniquely determined by the
requirement that $E=\bigoplus_{p\geq 0}\Lambda ^{(p,0)}(X)$ is $g^{\prime }-$%
orthogonal and $g^{\prime }=e^{f_{p}}g$ on $\Lambda ^{(p,0)}(X)$, where $%
f_{p}$ is a smooth function on $X$ to be determined. The decomposition $%
E=\bigoplus_{p\geq 0}\Lambda ^{(p,0)}(X)$ is $g^{\prime }-$orthogonal and
totally geodesic and we conclude that for $g^{\prime }$ we have $i\Lambda
\Theta _{g^{\prime }}\mid _{\Lambda ^{(p,0)}(X)}=i\Lambda \Theta _{(\Lambda
^{(p,0)}(X),g^{\prime }\mid \Lambda ^{(p,0)}(X))}=(\overline{\square }%
_{g}(f_{p})+\binom{n-1}{p-1}\lambda )Id_{\Lambda ^{(p,0)}(X)}$ on $\Lambda
^{(p,0)}(X)$ .

For the Higgs structure $\theta $ defined by the form $\varpi $ one checks
that $\overline{\theta }_{g}=\overline{\theta }_{g^{\prime }}$ because one
has just changed the metric conformally on each of the orthogonal subspaces
of $E$. Consequently we have $T_{g}=T_{g^{\prime }}$ for the operator
defined as in \ref{T-op}. Also because $E=\bigoplus_{p\geq 0}\Lambda
^{(p,0)}(X)$ is both $g$ and $g^{\prime }$ orthogonal, $T_{g}=T_{g^{\prime
}}:\Lambda ^{(p,0)}(X)\rightarrow \Lambda ^{(p,0)}(X)\forall p.$ Because $%
i\Lambda F_{g}=i\Lambda \Theta _{g}+T_{g}$, combining all these observations
yields 
\begin{equation}
i\Lambda F_{g^{\prime }}\mid _{\Lambda ^{(p,0)}(X)}=(\overline{\square }%
_{g}(f_{p})+\binom{n-1}{p-1}\lambda )Id_{\Lambda ^{(p,0)}(X)}+T_{g}\mid
_{\Lambda ^{(p,0)}(X)}  \label{HHYMex}
\end{equation}

Now using the $(n,0)-$form $\varpi $, which we can assume has pointwise
length $1$, we see from \ref{f[i] values} that $T_{g}\mid _{\Lambda
^{(p,0)}(X)}=$\bigskip $f[p]Id_{\Lambda ^{(p,0)}(X)}$ in the notation of \ref
{f[i] values}. Because of assumption ii above, we can find, for each $p$,
and for any constant $C,$a function $f_{p}$ such that 
\begin{equation}
\overline{\square }_{g}(f_{p})=C-\binom{n-1}{p-1}\lambda -f[p]  \label{HHYMf}
\end{equation}
Hence with such a choice of $f_{p}$ for each $p$, the corresponding
Hermitian metric $g^{\prime }$ on $E$ is Higgs-Hermitian-Yang-Mills, with
constant $C.$ Note that in general $g^{\prime }$ (restricted to $\Lambda
^{(1,0)}(X)$ and then defined on $T^{(1,0)}(X)$ by $g^{\prime }-$duality),
will not be a K\"{a}hler metric.

We now return to the question of the existence of HHYM metrics in the case
where $X$ a compact K\"{a}hler manifold. We do not have any examples of such
metrics for the Higgs bundles defined by \ref{higgsform2}. In fact, the
results we prove below on the nonexistence of such metrics came about as
obstructions to such metrics in our investigations of this question . We
need formulas for $c_{1}(F)$ for various subbundles of $E$. If X is any
complex manifold of complex dimension n. Then there is the well-known
formula 
\begin{equation}
c_{1}(\Lambda ^{(p,0)}(X))=\binom{n-1}{p-1}c_{1}(\Lambda ^{(1,0)}(X))
\label{c1}
\end{equation}

If $p=0$, we interpret $\binom{n-1}{p-1}$ to mean $0$, so the formula is
correct in this case, too. If we now assume $X$ is a compact K\"{a}hler
manifold, then it follows from \ref{c1} that $\deg (\Lambda ^{(p,0)}(X))=%
\binom{n-1}{p-1}\deg (\Lambda ^{(1,0)}(X))$ and $\mu (\Lambda
^{(p,0)}(X))=p\mu (\Lambda ^{(1,0)}(X))$.

Considering $E=\bigoplus_{p=0}^{n}\Lambda ^{(p,0)}(X)$ as a Higgs bundle via 
\ref{higgsform}, there are a large number of Higgs subbundles, hence Higgs
subsheaves of $\vartheta E$. For example, $\bigoplus_{p\geq 0}\Lambda
^{(2p,0)}(X)=:\Lambda ^{even}$ and $\bigoplus_{p\geq 0}\Lambda
^{(2p+1,0)}(X)=:\Lambda ^{odd}$ are both Higgs subbundles of $E$. Also,
using the fact that the first Chern class is additive over direct sums of
bundles, one computes that $c_{1}(E)=2^{n-1}c_{1}(\Lambda ^{(1,0)})$. One
computes $\mu (E)=\frac{n}{2}\mu (\Lambda ^{(1,0)})$, $c_{1}(\Lambda
^{even})=2^{n-2}c_{1}(\Lambda ^{(1,0)})=c_{1}(\Lambda ^{odd})$, and finally $%
\mu (\Lambda ^{even})=\frac{n}{2}\mu (\Lambda ^{(1,0)})=\mu (\Lambda
^{odd})=\mu (E)$. One gets a Higgs ``filtration'' $\{E^{a}\}_{a=0}^{n}$(i.e.
a filtration by Higgs subbundles) of $E$ as follows: 
\begin{gather}
0\subset E^{n}\subset E^{n-1}\subset \cdots \subset E^{1}\subset E^{0}=E
\label{filt1} \\
\text{where }E^{a}=\bigoplus_{p\geq a}\Lambda ^{(p,0)}(X)  \notag
\end{gather}
This filtration also gives Higgs filtrations of $\Lambda ^{even(odd)}$ by $%
\Lambda ^{even(odd)}\cap E^{a}$. Writing $k=2b+1$, each of the subbundles 
\begin{equation}
\bigoplus_{p\geq 0}\Lambda ^{(2bp+i,0)}(X),i=1,2,...,2b-1  \label{splitting}
\end{equation}
are Higgs subbundles of $\Lambda ^{even}$ or $\Lambda ^{odd}$, and
intersecting with the Higgs filtration gives more Higgs subbundles.

We now investigate the of \ the stability of some of these Higgs bundles. If 
$V\rightarrow X$ \ is a stable holomorphic vector bundle, then $V$ cannot
split holomorphically and nontrivially ($V=V_{1}\bigoplus V_{2}$ holomorphic
implies one $V_{i}=0)$, i.e. $V$ is irreducible$.$

This irreducibility result also holds for Higgs bundles: if $E$ is Higgs
stable, and $E=E_{1}\bigoplus E_{2}$ with $E_{1}$ and $E_{2}$ Higgs, then
one of these subbundles is $0$ (and the proof follows that of Lemma (7.3)
chapter 5 section 7 in \cite{kobayashi}). As a result, for the $E$ defined
by \ref{higgsform}, if $F$ is a Higgs subbundle of $E$ for which there is a
nontrivial splitting of the form $F=F\cap \Lambda ^{even}\bigoplus F\cap
\Lambda ^{odd}$ , then $F$ cannot be Higgs stable (or ``plain'' stable).
This type of splitting occurs in the bundles in the Higgs filtration \ref
{filt1}. In particular, none of the $E^{a},a=0,1,...,n-1$ can be stable,
although $E^{n},$ being a line bundle, is stable (cf. Proposition (7.7),
page 170 of \cite{kobayashi}).

It is natural to ask if any of the these Higgs subbundles could be
semistable. If any of the components of the Higgs filtration were
semistable(stable), say $E^{a},$then $\mu (E^{b})\leq (<)\mu (E^{a})$ $%
\forall b>a$. The next result shows that if $\deg (\Lambda ^{(1,0)}(X))>0$,
then the bundles $E^{a}$, $E^{a}\cap \Lambda ^{even}$, $E^{a}\cap \Lambda
^{odd}$ and others cannot be semistable in the ''ordinary'' sense where no
Higgs structure is assumed (again excluding the automatic case $E^{n}$ ,
which is a line bundle and hence stable), and will be used to show that many
of these bundles cannot admit HHYM metrics.

\begin{proposition}
\label{Eastable}Let $d=\deg (\Lambda ^{(1,0)}(X)),$ let $P$ be the
holomorphic subbundle of $\bigoplus_{i=0}^{n}\Lambda ^{(i,0)}(X)$ given by $%
P=\bigoplus_{s=1}^{z}\Lambda ^{(p_{s},0)}(X)$, $0\leq p_{1}<p_{2}<\cdots
<p_{z}\leq n$, and let $Q$ be the holomorphic subbundle of $P$ given by $%
Q=\bigoplus_{t=1}^{l}\Lambda ^{(q_{t},0)}(X)$, $p_{1}\leq q_{1}<q_{2}<\cdots
<q_{l}\leq p_{z}$, $\{q_{1},q_{2},...,q_{l}\}\subset
\{p_{1},p_{2},...,p_{z}\}$. If $Q$ is the ''tail'' of $P$, $%
q_{i}=p_{z-l+i},i=1,...,l$, then $\mu (P)>\mu (Q)\Leftrightarrow d<0,\mu
(P)=\mu (Q)\Leftrightarrow d=0$
\end{proposition}

\begin{proof}
We will show that $\mu (P)-\mu (Q)=d$ $c(p_{1},...,p_{z};q_{1},...q_{l})$
where\newline
$c(p_{1},...,p_{z};q_{1},...q_{l})$ is a rational number which is strictly
negative when $q_{i}=p_{z-l+i},i=1,...,l$. Write $\{p_{1},p_{2},...,p_{z}\}=%
\{\{q_{1},q_{2},...,q_{l}\},\{r_{1},r_{2},...,r_{z-l}\}\},r_{1}<r_{2}<\cdots
<r_{z-l}$. Using the formula \ref{c1} one computes 
\begin{gather}
\mu (P)=\frac{\sum_{s=1}^{z}\binom{n-1}{p_{s}-1}d}{\sum_{s=1}^{z}\binom{n}{%
p_{s}}}  \notag \\
\mu (P)-\mu (Q)=\frac{d}{rk(P)rk(Q)}\left( \sum_{s=1}^{z}\binom{n-1}{p_{s}-1}%
\sum_{t=1}^{l}\binom{n}{q_{t}}-\sum_{t=1}^{l}\binom{n-1}{q_{t}-1}%
\sum_{s=1}^{z}\binom{n}{p_{s}}\right)  \notag \\
=\frac{d}{rk(P)rk(Q)}\left( \sum_{b=1}^{z-l}\binom{n-1}{r_{b}-1}%
\sum_{t=1}^{l}\binom{n}{q_{t}}-\sum_{t=1}^{l}\binom{n-1}{q_{t}-1}%
\sum_{b=1}^{z-l}\binom{n}{r_{b}}\right)  \notag \\
=\frac{d}{rk(P)rk(Q)}\sum_{\substack{ 1\leq t\leq l  \\ 1\leq b\leq z-l}}%
\binom{n-1}{r_{b}-1}\binom{n-1}{q_{t}-1}n\{\frac{r_{b}-q_{t}}{r_{b}\,q_{t}}\}
\label{diff}
\end{gather}

In the last line of \ref{diff} we have assumed $r_{1},\,q_{1}\geq 1$. If $Q$
is the tail of $P$, then $r_{b}-q_{t}<0$ for all $b$ and $t$. In case one or
both of $r_{1}$ or $\,q_{1}$ is $0$, one has to write out some special cases
of the expression in \ref{diff}, but the basic result is the same: $\mu
(P)-\mu (Q)=d$ $c$ where $c$ is a rational number, which is strictly
negative if $r_{b}-q_{t}<0$ for all $b$ and $t.$
\end{proof}

\begin{theorem}
\label{bigthm}Let $X$ be a compact K\"{a}hler manifold with a nontrivial
holomorphic k-form $\varpi $ where $k>1$ is odd. Let the Higgs structure of $%
E$ be as above. and let $P$ be any Higgs subbundle of $E$ of the form $%
P=\bigoplus_{s=1}^{z}\Lambda ^{(p_{s},0)}(X)$, $0\leq p_{1}<p_{2}<\cdots
<p_{z}\leq n$, $(z\geq 2)$. Then $P$ does not admit any
Higgs-Hermitian-Yang-Mills metric in any of the following cases :
\end{theorem}

i. $\deg (X)<0$

ii. $\deg (X)=0$ and $p_{1}\leq n-k+1$

iii. $k\geq \frac{n}{2}+1$, $p_{1}\leq n-k+1$, and $\varpi $ is a section of 
$P$.

\begin{proof}
We first prove i. \ Note that $\deg (\Lambda ^{(s,0)}(X))=\binom{n-1}{s-1}%
\deg (\Lambda ^{(1,0)}(X))$, hence $\deg (P)$ is a positive multiple of $d=$ 
$\deg (\Lambda ^{(1,0)}(X))=-\deg (X)$. Assume $P$ admits a HHYM metric $h,$ 
$i\Lambda F_{h}=cId_{P}.$ \ Because $P$ admits this HHYM metric, it is
Higgs-semistable, so $\mu (P)\geqslant $ $\mu (P^{\prime })$, where $%
P^{\prime }$ is the Higgs subbundle of $P$ given by $\bigoplus_{s=2}^{z}%
\Lambda ^{(p_{s},0)}(X)$. \ Now using $P,$ and $Q=P^{\prime }$ in
Proposition \ref{Eastable}, we get that $-\deg (X)=d\leq 0$, proving i. We
now prove ii. \ Assume that $P$ admits a HHYM metric$\ h,$ $i\Lambda F_{h}$ $%
=cId_{P}$. Representing $c_{1}(P)$ by $\frac{i}{2\pi }tr_{P}F_{h}=\frac{i}{%
2\pi }\sum_{\alpha =1}^{rkP}F_{h\alpha }^{\alpha }$ one computes $%
c_{1}(P)\wedge \omega ^{n-1}=\frac{1}{2\pi n}\sum_{\alpha }i\Lambda
F_{h\alpha }^{\alpha }\omega ^{n}$ $=\frac{rk(P)c}{2\pi n}\omega ^{n}$(\cite
{kobayashi} chapter 3, section 1, (1.18)) hence $\deg (P)$ is a positive
multiple of $c$.

If $d=0$, then we have $\mu (P)=0$ and also $\mu (P^{\prime })=0.$ If one
adapts the proof of Proposition (8.2), Chapter V of \cite{kobayashi} to the
Higgs setting one concludes the following, using the notation in \cite
{kobayashi}: if $E$ is a Higgs bundle and $E^{\prime }\subset E$ a Higgs
subbundle, over a compact K\"{a}hler manifold $(X,g),$ and if $E$ admits a
HHYM metric $h$, then $\mu (E^{\prime })\leq $ $\mu (E)$, with equality iff $%
E=E^{\prime }\oplus (E^{\prime })^{\bot _{h}}$ is a holomorphic splitting
into Higgs subbundles (i.e. $(E^{\prime })^{\bot _{h}}$ is a holomorphic,
Higgs subbundle of $E$ ). In our setting ($d=0,$so $\mu (P)=$ $\mu
(P^{\prime })=0$) this fact implies that $P=P^{\prime }\oplus (P^{\prime
})^{\bot _{h}}$ is a holomorphic Higgs splitting, so $(P^{\prime })^{\bot
_{h}}$ is $\theta $-invariant. However, $\theta P\subset P^{\prime }\otimes
\Lambda ^{(1,0)}(X)$. Therefore $\theta $ $(P^{\prime })^{\bot _{h}}=0$. Now
let $s\in \Gamma \Lambda ^{(p_{1},0)}(X)$ and split $s=s^{\prime }+s^{\prime
\prime }$, $s^{\prime }\in \Gamma P^{\prime }$, $s^{\prime \prime }\in
\Gamma (P^{\prime })^{\bot _{h}}$. Then $\theta s\in \Gamma \Lambda
^{(p_{1}+k-1,0)}(X)\otimes \Lambda ^{(1,0)}(X)$, while $\theta (s^{\prime
}+s^{\prime \prime })=\theta s^{\prime }\in \Gamma P^{\prime }\otimes
\Lambda ^{(1,0)}(X)$, so the lowest possible degree ''form'' part of $\theta
s^{\prime }$ is $p_{2}+k-1>p_{1}+k-1$. \ We conclude that $\theta s=0$, for
every $s\in \Gamma \Lambda ^{(p_{1},0)}(X)$. From Proposition \ref{kertheta}%
, we conclude that $\varpi \wedge i(Z)s=0$ for every $s\in \Gamma \Lambda
^{(p_{1},0)}(X)$ and every holomorphic tangent vector $Z$. Because $%
p_{1}\leq n-k+1$, this implies that $\varpi \equiv 0$ (one can see this
pointwise by picking $s=dz_{A}$, where $A=\{1\leq A_{1}<\cdots
<A_{p_{1}}\leq n\}$ can be any length $p_{1}$ multiindex, then $%
i(Z_{A_{p_{1}}})s=\pm dz_{A_{1}}\wedge dz_{A_{2}}\wedge dz_{A_{3}}\wedge
\cdots dz_{A_{p_{1}-1}}$ so any of the simple $p_{1}-1$ forms $%
dz_{A_{1}}\wedge dz_{A_{2}}\wedge dz_{A_{3}}\wedge \cdots dz_{A_{p_{1}-1}}$
can be obtained as $i(Z)s$). We have reached a contradiction.

Now assume $k\geq \frac{n}{2}+1$, but not necessarily that $\deg (X)=0$, and
as above, assume $P$ admits a HHYM metric $h,$ $i\Lambda F_{h}=cId_{P}.$
Because $P$ admits this HHYM metric, it is semistable and using $P,$ and $%
Q=P^{\prime }$ in Proposition \ref{Eastable} we get, that $d\leq 0$ so $%
c\leq 0.$ Now $c\leq 0$ implies that $i\Lambda F_{h}$ is a pointwise
nonpositive operator. Since $\varpi $ is a section of $P,$ the formula \ref
{boch5} in the current setting, for the bundle $P$, with $s=\varpi $, (cf.
Remark \ref{rmksbbdl}) becomes 
\begin{equation*}
\overset{}{c\Vert \varpi \Vert _{h}^{2}=-i\Lambda \partial \overline{%
\partial }\Vert \varpi \Vert _{h}^{2}+\Vert \nabla \varpi \Vert
_{h}^{2}+\sum \{}\Vert \overline{\theta _{h}}(\overline{Z_{i}})\varpi \Vert
_{h}^{2}\}
\end{equation*}
because $k\geq \frac{n}{2}+1$ implies $\theta \varpi =0$. Integrating this
equality over $X$ implies $c\geq 0$, hence $c=0$ and $\deg (X)=0.$ Now part
ii. gives a contradiction.
\end{proof}

\section{Section 3}

\refstepcounter{section}

In this section we examine analogs of the Kodaira and Nakano-type vanishing
theorems for $(p,q)$-forms with values in a Higgs bundle over a K\"{a}hler
manifold, cf. \cite{kobayashi} chapter 3 (see also \cite{shiffman} chapter
1) for the K\"{a}hler manifold operators and \cite{simpIHES} section 1 for
the formulas on Higgs bundles.

Let $(X,g)$ be a complex manifold of complex dimension $n$ with a K\"{a}hler
metric $g$. Let $(E,h)\rightarrow X$ be a rank r holomorphic vector bundle
with a Hermitian metric $h$. As in \cite{kobayashi}, the Hermitian
connection on $E$ extends to an operator $d^{\nabla }:C^{\infty }\Gamma
E\bigotimes \Lambda ^{i}(X)\rightarrow C^{\infty }\Gamma E\bigotimes \Lambda
^{i+1}(X)$ and there is the refinement of $d^{\nabla }$ into the two
operators $d^{\nabla }=\partial _{h}+\overline{\partial }$%
\begin{gather*}
\partial _{h}:C^{\infty }\Gamma E\bigotimes \Lambda ^{(i,j)}(X)\rightarrow
C^{\infty }\Gamma E\bigotimes \Lambda ^{(i+1,j)}(X) \\
\overline{\partial }:C^{\infty }\Gamma E\bigotimes \Lambda
^{(i,j)}(X)\rightarrow C^{\infty }\Gamma E\bigotimes \Lambda ^{(i,j+1)}(X)
\end{gather*}
given relative to a local holomorphic frame $\{e_{\alpha }\}_{\alpha =1}^{r}$
of $E$ by 
\begin{gather*}
\overline{\partial }(e_{\alpha }\bigotimes \phi )=e_{\alpha }\bigotimes 
\overline{\partial }\phi \\
\partial _{h}(e_{\alpha }\bigotimes \phi )=\overset{r}{\underset{\beta =1}{%
\sum }}e_{\beta }\bigotimes C_{\alpha }^{\beta }\wedge \phi +e_{\alpha
}\bigotimes \partial \phi
\end{gather*}
where $\nabla e_{\alpha }=\overset{r}{\underset{\beta =1}{\sum }}$ $e_{\beta
}\bigotimes C_{\alpha }^{\beta }$. The metric on $X$ extends to $\Lambda
^{\ast }(X)\bigotimes \mathbb{C}$ (we drop the $\mathbb{C}$ and write $%
\Lambda ^{\ast }(X)$ for the complex exterior algebra of $X$) and the
Hermitian metric on $E$ combines to give a metric which we denote $\ll $ , $%
\gg $ on $E\bigotimes \Lambda ^{\ast }(X)$ by the prescription 
\begin{gather*}
\ll e\bigotimes \phi ,f\bigotimes \psi \gg =h(e,f)\ast (\phi \wedge \ast 
\overline{\psi }) \\
e,f\in E\text{ }\phi ,\psi \in \Lambda ^{\ast }(X)
\end{gather*}
where $\ast $ denotes the Hodge star operator on $\Lambda ^{\ast }(X)$
determined by $g$ (and extend $\ast $ to $E\bigotimes \Lambda ^{\ast }(X)$
by $id_{E}\bigotimes \ast $ which we also denote simply as $\ast )$. We use
the convention that the Hodge star operator is given on the complex exterior
algebra by the method in \cite{kobayashi} Chapter 3, section 2. Note that
there is a sign typographical error in this reference in formula (2.6),
which should read 
\begin{equation*}
\varepsilon (A,B)=(-1)^{np+n(n-1)/2}\sigma (AA^{\prime })\cdot \sigma
(BB^{\prime })
\end{equation*}
The $L^{2}$ or formal adjoints of $\partial ^{\nabla }$ and $\overline{%
\partial }$ with respect to $\ll $ , $\gg $ are given by (cf \cite{kobayashi}
chapter 3 section 2) 
\begin{gather*}
\partial _{h}^{\ast }=-\ast \overline{\partial }\ast =i[\Lambda ,\overline{%
\partial }] \\
\overline{\partial }^{\ast }=-\ast \partial _{h}\ast =-i[\Lambda ,\partial
_{h}]
\end{gather*}
where $\Lambda =i\sum_{i,j}g^{i\overline{j}}i(\frac{\partial }{\partial z_{i}%
})i(\frac{\partial }{\partial \overline{z_{j}}})$ is the adjoint to exterior
multiplication with the K\"{a}hler form $\omega $. Let $\square _{\overline{%
\partial }}$ $=(\overline{\partial }+\overline{\partial }^{\ast })^{2}$ and $%
\square _{\partial _{h}}=(\partial _{h}+\partial _{h}^{\ast })^{2}$. The
Kodaira-Nakano formula (cf. \cite{kobayashi} see Chapter 3 section 3, the
proof of (3.5) page 69, or \cite{shiffman}, Chapter 1, page 16, (1.58)) can
be written 
\begin{equation}
\square _{\partial _{h}}-\square _{\overline{\partial }}=i(\Lambda e(\Theta
)-e(\Theta )\Lambda )  \label{nakano}
\end{equation}
where as in \cite{kobayashi} $e(\Theta ):E\bigotimes \Lambda
^{(p,q)}(X)\rightarrow E\bigotimes \Lambda ^{(p+1,q+1)}(X)$ is given by $%
e(\Theta )(e_{\alpha }\bigotimes \phi )=\underset{\beta }{\sum }\{e_{\beta
}\otimes \Theta _{\alpha }^{\beta }\wedge \phi \}$. We will say that a
section $s$ of $E\bigotimes \Lambda ^{w}(X)=\bigoplus_{i+j=w}E\bigotimes
\Lambda ^{(i,j)}(X)$, has Hodge type $(p,q)$ if $s=\sum_{\alpha
=1}^{r}e_{\alpha }\bigotimes \phi _{\alpha }$ where each $\phi _{\alpha }$
is a section of $\Lambda ^{(p,q)}(X)$. This terminology can become ambiguous
if $E$ is the Higgs bundle discussed in section 2, since the $E$ component
of a section of$\ E\bigotimes \Lambda ^{(p,q)}(X)$ will be a sum of forms
with Hodge types. We will address this issue when it arises. Both $\square
_{\partial _{h}}$ and $\square _{\overline{\partial }}$\ preserve the Hodge $%
(p,q)$ types.

Now suppose $E$ also has the structure of a Higgs bundle with Higgs form $%
\theta $. The operators $D^{\prime \prime }=\overline{\partial }+\theta
,D_{h}^{\prime }$ and $D_{h}$ defined in \ref{D'h} and \ref{Dh} extend to $%
E\bigotimes \Lambda ^{\ast }(X)$ as above (\cite{simpIHES} section 1):

\begin{gather}
D^{\prime \prime }(e_{\alpha }\bigotimes \phi )=e_{\alpha }\bigotimes 
\overline{\partial }\phi +\sum_{\beta }e_{\beta }\bigotimes \theta _{\alpha
}^{\beta }\wedge \phi  \label{Dext} \\
D_{h}^{\prime }(e_{\alpha }\bigotimes \phi )=e_{\alpha }\bigotimes \partial
\phi +\sum_{\beta }e_{\beta }\bigotimes (C_{\alpha }^{\beta }+\overline{%
\theta _{h}}_{\alpha }^{\beta })\wedge \phi  \notag \\
D_{h}^{\prime \ast }=-\ast D^{\prime \prime }\ast =i[\Lambda ,D^{\prime
\prime }]  \notag \\
D^{\prime \prime \ast }=-\ast D_{h}^{\prime }\ast =-i[\Lambda
,D_{h}^{^{\prime }}]  \notag \\
D_{h}^{\ast }=D_{h}^{\prime \ast }+D^{\prime \prime \ast }  \notag
\end{gather}

One checks that as before that these extended operators $D^{\prime \prime
},D_{h}^{\prime }$ and their adjoints all square to zero. Note the adjoints
of the Higgs forms are given by $\theta ^{\ast }=\ast \overline{\theta }%
_{h}\ast =-i[\Lambda ,\overline{\theta }_{h}]$, $\overline{\theta }%
_{h}^{\ast }=\ast \theta \ast =i[\Lambda ,\theta ]$. Define 
\begin{gather}
\square _{D^{\prime \prime }}=(D^{\prime \prime }+D^{\prime \prime \ast
})^{2}=D^{\prime \prime }D^{\prime \prime \ast }+D^{\prime \prime \ast
}D^{\prime \prime }  \label{D''laplacian} \\
\square _{D_{h}^{^{\prime }}}=(D_{h}^{\prime }+D_{h}^{\prime \ast
})^{2}=D_{h}^{\prime }D_{h}^{\prime \ast }+D_{h}^{\prime \ast }D_{h}^{\prime
}  \label{D'laplacian} \\
\square _{D_{h}}=D_{h}D_{h}^{\ast }+D_{h}^{\ast }D_{h}=\square _{D^{\prime
\prime }}+\square _{D_{h}^{^{\prime }}}  \notag
\end{gather}

We now examine the relation between the Laplacians in \ref{D'laplacian} and 
\ref{D''laplacian}, and the ``ordinary'' Laplacians corresponding to $\theta
=0$, $\square _{\partial _{h}},\square _{\overline{\partial }}$ acting on $%
E\bigotimes \Lambda ^{t}(X)$. The following formulas are given by
computations using the definitions of $D^{\prime \prime }$, $D_{h}^{^{\prime
}}$ and their adjoints \ref{Dext} (cf. the notation discussed at the
beginning of the proof of Proposition \ref{prop}) : 
\begin{eqnarray}
\square _{D^{\prime \prime }} &=&\square _{\overline{\partial }}+\theta
\theta ^{\ast }+\theta ^{\ast }\theta +\overline{\partial }(\theta ^{\ast })+%
\overline{\partial }^{\ast }(\theta )  \label{lap} \\
\square _{D_{h}^{^{\prime }}} &=&\square _{\partial _{h}}+\overline{\theta }%
_{h}\overline{\theta }_{h}^{\ast }+\overline{\theta }_{h}^{\ast }\overline{%
\theta }_{h}+\partial _{h}(\overline{\theta }_{h}^{\ast })+\partial
_{h}^{\ast }(\overline{\theta }_{h})  \notag
\end{eqnarray}

In general the operators $\square _{D^{\prime \prime }}$ and $\square
_{D_{h}^{^{\prime }}}$ will not preserve the Hodge type $(p,q)$ of a section
of $E\bigotimes \Lambda ^{(p,q)}(X)$, although they do preserve the total
degree $p+q$ because 
\begin{gather}
\overline{\partial }(\theta ^{\ast }),\partial _{h}^{\ast }(\overline{\theta 
}_{h}):C^{\infty }\Gamma E\bigotimes \Lambda ^{(p,q)}(X)\rightarrow
C^{\infty }\Gamma E\bigotimes \Lambda ^{(p-1,q+1)}(X)  \label{zerothord} \\
\overline{\partial }^{\ast }(\theta ),\partial _{h}(\overline{\theta }%
_{h}^{\ast }):C^{\infty }\Gamma E\bigotimes \Lambda ^{(p,q)}(X)\rightarrow
C^{\infty }\Gamma E\bigotimes \Lambda ^{(p+1,q-1)}(X)  \notag
\end{gather}
and defining 
\begin{gather*}
\overline{\boxminus }\text{ := }\square _{\overline{\partial }}+\theta
\theta ^{\ast }+\theta ^{\ast }\theta \\
\boxminus \text{:= }\square _{\partial _{h}}+\overline{\theta }_{h}\overline{%
\theta }_{h}^{\ast }+\overline{\theta }_{h}^{\ast }\overline{\theta }_{h}
\end{gather*}
each of $\boxminus $ and $\overline{\boxminus }$ preserve Hodge type$(p,q)$
of a section of $E\bigotimes \Lambda ^{w}(X)$ ,i.e. $\boxminus $, $\overline{%
\boxminus }:$ $E\bigotimes \Lambda ^{(p,q)}(X)\rightarrow E\bigotimes
\Lambda ^{(p,q)}.$ The operators $\boxminus $, and $\overline{\boxminus }$
differ from the usual ($\theta =0$) Laplacians only by the zeroth order
terms, which are nonnegative operators$.$ If $X$ is compact K\"{a}hler, then
because $\overline{\boxminus }$ and $\boxminus $ preserve Hodge type $(p,q)$
one has $\ker \overline{\boxminus }(\boxminus )\subset \ker \square
_{D^{\prime \prime }}(\square _{D_{h}^{^{\prime }}})$ considering these
operators acting on $C^{\infty }\Gamma E\bigotimes \Lambda ^{w}(X).$ To wit,
if $s=\sum_{p+q=w}s_{p,q}$ is the decomposition into Hodge $(p,q)$
components, then $\overline{\boxminus }s=0\Leftrightarrow \overline{%
\boxminus }s_{p,q}=0\forall (p,q)\Leftrightarrow 0=\overline{\partial }%
s_{p,q}=\theta s_{p,q}=\overline{\partial }^{\ast }s_{p,q}=\theta ^{\ast
}s_{p,q}\forall (p,q)\Leftrightarrow $\newline
$0=D^{\prime \prime }s_{p,q}=D^{\prime \prime \ast }s_{p,q}\forall
(p,q)\Rightarrow 0=D^{\prime \prime }s=D^{\prime \prime \ast
}s\Leftrightarrow \square _{D^{\prime \prime }}s.$

An analog of the Kodaira-Nakano-type formula in this setting is 
\begin{equation}
\square _{D_{h}^{^{\prime }}}-\square _{D^{\prime \prime }}=i(\Lambda
e(F_{h})-e(F_{h})\Lambda )  \label{higgsnak1}
\end{equation}
where as in \cite{kobayashi} $e(F_{h}):E\bigotimes \Lambda
^{w}(X)\rightarrow E\bigotimes \Lambda ^{w+2}(X)$ is given $%
e(F_{h})(e_{\alpha }\bigotimes \phi )=\underset{\beta }{\sum }\{e_{\beta
}\otimes F_{h\alpha }^{\beta }\wedge \phi \}$. We note that 
\begin{gather}
i(\Lambda e(F_{h}^{(1,1)})-e(F_{h}^{\left( 1,1\right) })\Lambda
):E\bigotimes \Lambda ^{(p,q)}(X)\rightarrow E\bigotimes \Lambda ^{(p,q)}
\label{kodnakcrv} \\
i(\Lambda e(F_{h}^{(2,0)})-e(F_{h}^{(2,0)})\Lambda ):E\bigotimes \Lambda
^{(p,q)}(X)\rightarrow E\bigotimes \Lambda ^{(p+1,q-1)}  \notag \\
i(\Lambda e(F_{h}^{(0,2)})-e(F_{h}^{(0,2)})\Lambda ):E\bigotimes \Lambda
^{(p,q)}(X)\rightarrow E\bigotimes \Lambda ^{(p-1,q+1)}  \notag
\end{gather}
Because $\square _{D^{\prime \prime }}$, $\square _{D_{h}^{^{\prime }}}$ and 
$\square _{D}$ are nonnegative operators on a compact K\"{a}hler manifold,
formula \ref{higgsnak1} implies

\begin{theorem}
\label{thm1}If ($X,g)$ is a compact K\"{a}hler manifold, $(E,\theta
)\rightarrow X$ a Higgs bundle with a Hermitian metric $h$ then with the
notation above 
\begin{gather}
i(\Lambda e(F_{h})-e(F_{h})\Lambda )\leq 0\Rightarrow  \notag
\label{vanish2} \\
\ker \square _{D^{\prime \prime }}\subseteq \ker \square _{D_{h}^{^{\prime
}}}\cap \ker \square _{D}
\end{gather}
where $i(\Lambda e(F_{h})-e(F_{h})\Lambda ),\square _{D^{\prime \prime
}},\square _{D_{h}^{^{\prime }}}$ and $\square _{D}$ are considered as
operators on $C^{\infty }\Gamma E\bigotimes \Lambda ^{w}(X).$
\end{theorem}

\begin{proof}
That $\ker \square _{D^{\prime \prime }}\subseteq \ker \square
_{D_{h}^{^{\prime }}}$ follows from \ref{higgsnak1}. Then, if $s\in \ker
\square _{D^{\prime \prime }}\cap \ker \square _{D_{h}^{^{\prime }}}$, one
gets $0=D^{\prime \prime }s=D^{\prime \prime \ast }s=D_{h}^{\prime
}s=D_{h}^{\prime \ast }s$ and hence $0=Ds=D^{\ast }s$.
\end{proof}

If the Higgs operator $\theta $ is parallel with respect to the operator $%
d^{\nabla }$ defined by $h$ as in Proposition \ref{prop}, as a section of $%
Hom(E)\bigotimes \Lambda ^{(1,0)}(X),$ i.e. $F_{h}$ has only $(1,1)$ form
parts relative to a holomorphic frame, then the four operators in \ref
{zerothord} all vanish (the proof is analogous to the computation in
Proposition \ref{prop}) and $\square _{D^{\prime \prime }}=\overline{%
\boxminus }$ and $\square _{D_{h}^{\prime }}=\boxminus $ preserve Hodge type.

Let $s=\sum_{p+q=w}s_{p,q}\in C^{\infty }\Gamma E\bigotimes \Lambda ^{w}(X)$
be as above, then the formulas above combine to give 
\begin{gather}
\ll (\square _{D_{h}^{^{\prime }}}-\square _{D^{\prime \prime }})s,s\gg =\ll
(i(\Lambda e(F_{h})-e(F_{h})\Lambda )s,s\gg =  \notag \\
\sum_{p+q=w}\{\ll (\boxminus -\overline{\boxminus })s_{p,q},s_{p,q}\gg
\label{genlap} \\
\ll (\partial _{h}(\overline{\theta }_{h}^{\ast })-\overline{\partial }%
^{\ast }(\theta )_{h})s_{p,q},s_{p+1,q-1}\gg +\ll (\partial _{h}^{\ast }(%
\overline{\theta }_{h})-\overline{\partial }(\theta ^{\ast
}))s_{p,q},s_{p-1,q+1}\gg \}  \notag
\end{gather}

In formula \ref{genlap} we see that if $s\in E\bigotimes \Lambda ^{(p,q)}(X)$
then the term $\ll (i(\Lambda e(F_{h})-e(F_{h})\Lambda )s,s\gg $ depends
only on the $(1,1)$ part of $F_{h}$ cf. \ref{kodnakcrv}. This observation
yields the following vanishing result, which will be revisited in giving a
type of Higgs bundle analog to the Kodaira-Nakano vanishing theorem.

\begin{theorem}
\label{thm3.2}Let ($X,g)$ be a compact K\"{a}hler manifold, $(E,\theta
)\rightarrow X$ a Higgs bundle with a Hermitian metric $h$. With the
notation above, if the $(1,1)$ part of $i(\Lambda e(F_{h})-e(F_{h})\Lambda
)\leq 0$ pointwise as an operator on $E\bigotimes \Lambda ^{(p,q)}(X)$ and
if $s_{p,q}\in C^{\infty }\Gamma E\bigotimes \Lambda ^{(p,q)}(X)$ satisfies $%
\square _{D^{\prime \prime }}s_{p,q}=0$ then $\square _{D_{h}^{^{\prime
}}}s_{p,q}=0$ and $\square _{D_{h}}s_{p,q}=0$. If the $(1,1)-$part of $%
i(\Lambda e(F_{h})-e(F_{h})\Lambda )$ is quasinegative on $E\bigotimes
\Lambda ^{(p,q)}(X)$, any then any such $s_{p,q}$ must be $0.$\newline
\qquad\ If the Higgs form $\theta $ is parallel as a section of $%
Hom(E)\bigotimes \Lambda ^{(1,0)}(X),$i.e. $F_{h}$ has only $(1,1)$ form
parts with respect to a holomorphic frame, and if $i(\Lambda
e(F_{h})-e(F_{h})\Lambda )\leq 0$ on $E\bigotimes \Lambda ^{w}(X)$ then $%
s=\sum_{p+q=w}s_{p,q}\in E\bigotimes \Lambda ^{w}(X),\square _{D^{\prime
\prime }}s=0$ implies $\square _{D_{h}^{^{\prime }}}s_{p,q}=0$ and $\square
_{D_{h}}s_{p,q}=0$ $\forall (p,q)$(cf. Theorem \ref{thm1}).
\end{theorem}

The Kodaira-Nakano vanishing theorems are generally stated as vanishing
theorems for harmonic sections of $(p,q)$-forms with values in a holomorphic
line bundle $L$ i.e. harmonic sections of $L\bigotimes \Lambda ^{(p,q)}(X)$,
with $X$ compact K\"{a}hler or compact complex with $c_{1}(L)<0$ (\cite
{kobayashi} chapter 3 section 3, and \cite{shiffman} chapter 2 Theorem
(2.18)). However the proofs given work for $\square _{\overline{\partial }}-$%
harmonic sections of $(p,q)$-forms with values in a Hermitian holomorphic
vector bundle, i.e. harmonic sections of $E\bigotimes \Lambda ^{(p,q)}(X)$,
if we assume that $E$ admits a projectively flat Hermitian metric (every
Hermitian metric on a line bundle is projectively flat). The main technical
point is that conformally changing a projectively flat Hermitian metric
yields another projectively flat Hermitian metric. This observation is used
in the course of proving the next theorem.

\begin{theorem}
\label{kodnakseam}(Kodaira-Nakano) Let $(E,\theta )\rightarrow (X,g)$ be a
Higgs bundle over a compact complex manifold of complex dimension n with $%
c_{1}(E)<0$ and assume $E$ admits a Hermitian metric $h$ for which the $%
(1,1) $ part of the Higgs curvature $F_{h}^{(1,1)}$satisfies the equation $%
F_{h}^{(1,1)}=\kappa Id_{E}$, where $\kappa =\sum \kappa _{i\overline{j}%
}dz_{i}\wedge \overline{dz_{j}}$ is a $(1,1)$ form. If $s_{p,q}\in C^{\infty
}\Gamma E\bigotimes \Lambda ^{(p,q)}(X)$ and $\square _{D^{\prime \prime
}}s_{p,q}=0,$ $p+q\leq n-1(=n)$, then $s_{p,q}=0(\square _{D_{h}^{^{\prime
}}}s_{p,q}=0$ and $\square _{D_{h}}s_{p,q}=0)$.

\begin{proof}
Let a negative representative of $c_{1}(E)$ be given by $\frac{i}{2\pi }f=%
\frac{i}{2\pi }\sum_{i,j}f_{i\overline{j}}dz_{i}\wedge d\overline{z}_{j}$ ,
a closed real $(1,1)$ form with $(f_{i\overline{j}})$ negative definite
pointwise. Then \newline
$-i\sum_{i,j}f_{i\overline{j}}dz_{i}d\overline{z}_{j}$ is a K\"{a}hler
metric $g$ on $X$, and we use this metric as the K\"{a}hler metric on $X$.
Suppose $h$ is a Hermitian metric on $E$ for which $F_{h}^{(1,1)}=\kappa
Id_{E}$ where $\kappa =\sum \kappa _{i\overline{j}}dz_{i}\wedge \overline{%
dz_{j}}$ is a $\left( 1,1\right) $ form. From \ref{F11} we have $%
tr_{E}F_{h}^{(1,1)}=tr_{E}\Theta $ and thus we can represent $c_{1}(E)$ as$%
\frac{i}{2\pi }tr_{E}F_{h}^{(1,1)}=\frac{i}{2\pi }r$ $\kappa $. Now, any
conformal change of $h$\ to a new Hermitian metric $h^{\prime }$ $=ah$ ($a$
a smooth positive real-valued function on $X$) changes the Hermitian metric
curvature from $\Theta $ to $\Theta ^{\prime }=\Theta -\partial \overline{%
\partial }\ln (a)Id_{E}$ and does not change $\overline{\theta _{h}},$i.e. $%
\overline{\theta _{h^{\prime }}}$ $=\overline{\theta _{h}}.$ It therefore
follows from \ref{F11} that for the new metric $h^{\prime },F_{h^{\prime
}}^{(1,1)}=\kappa ^{\prime }Id_{E}$ where $\kappa ^{\prime }=$ $\kappa
-\partial \overline{\partial }\ln (a)$.

\qquad For any real representative of $c_{1}(E)$, such as $\frac{i}{2\pi }f,$
we can always change $h$ (or any given Hermitian metric on a $E)$
conformally to a new metric $h^{\prime }$ for which $tr_{E}\Theta ^{\prime
}=f$ (\cite{kobayashi} Chapter 2 section 2, page 41, Proposition (2.23)).
Therefore we conclude: we can conformally change $h$ to a metric $h^{\prime
} $ for which $f=tr_{E}\Theta ^{\prime }=tr_{E}F_{h^{\prime }}^{(1,1)}=r$ $%
\kappa ^{\prime }$. Thus $F_{h^{\prime }}^{(1,1)}=\frac{1}{r}fId_{E}$. Now
with this Hermitian metric the formulas \ref{higgsnak1} and \ref{kodnakcrv}
yield, if $s\in C^{\infty }\Gamma E\bigotimes \Lambda ^{(p,q)}(X),$ $%
h((\square _{D_{h}^{^{\prime }}}-\square _{D^{\prime \prime
}})s,s)=h(i(\Lambda e(F_{h})-e(F_{h})\Lambda )s,s)=h(i(\Lambda
e(F_{h}^{(1,1)})-e(F_{h}^{\left( 1,1\right) })\Lambda )s,s)=$\newline
$h(-\frac{1}{r}(\Lambda L-L\Lambda )s,s)=-\frac{1}{r}(n-(p+q))\left\|
s\right\| _{h^{\prime }}^{2}$ ($L$ is exterior multiplication by the
K\"{a}hler form). Thus if $p+q\leq n-1\left( =n\right) $, $\square
_{D^{\prime \prime }}s=0$ we conclude $s=0(\square _{D_{h}^{^{\prime }}}s=0$
and $\square _{D_{h}}s=0)$.
\end{proof}
\end{theorem}

\begin{remark}
The vanishing theorem of Gigante and Girbau, (\cite{kobayashi} chapter 3
section 3, Theorem (3.4) and \cite{shiffman} chapter 3 Theorem (3.2)) where
the assumptions are: $X$ is compact K\"{a}hler, $c_{1}(E)\leq 0$ and
pointwise rank $k,$ and the vanishing occurs in degrees $p+q\leq k-1$, is
also valid for a Hermitian holomorphic vector bundle with a projectively
flat metric $h$ and in the Higgs setting when $F_{h}^{(1,1)}=\kappa Id_{E}$.
The proof given in \cite{kobayashi} chapter 3 section 3, Theorem (3.4),
pages 69-73 works in this setting. One has to extend the formula (3.6), page
70 as we now indicate. Let $(E,\theta )\rightarrow (X,g)$ be a Higgs bundle
over a compact K\"{a}hler manifold and $E$ admits a Hermitian metric $h$ for
which $F_{h}^{(1,1)}=\kappa Id_{E}$, where $\kappa =\sum \kappa _{i\overline{%
j}}dz_{i}\wedge \overline{dz_{j}}$ is a $(1,1)$ form. Then $\kappa Id_{E}$
has the same Hermitian symmetries as the Hermitian curvature, $\Theta $, of $%
h$, so w.l.o.g., $\kappa =\sum \kappa _{i}dz_{i}\wedge \overline{dz_{i}}$ ($%
\{dz_{i}\}_{i=1}^{n}$ orthonormal at the point of evaluation). Now for $s\in
C^{\infty }\Gamma E\bigotimes \Lambda ^{(p,q)}(X)$ we can write locally $%
s=\sum e_{a}\bigotimes \varphi _{I\overline{J}}^{a}dz_{I}\wedge \overline{%
dz_{J}}$ where the multiindices satisfy $\left| I\right| =p,\left| J\right|
=q$ and at one particular point of evaluation the holomorphic frame $%
\{e_{a}\}_{a=1}^{r}$ is orthonormal. With this notation the formula (3.6),
page 70 in \cite{kobayashi} translates into 
\begin{equation*}
h(i(\Lambda e(F_{h}^{(1,1)})-e(F_{h}^{\left( 1,1\right) })\Lambda )s,s)=\sum 
_{\substack{ a=1  \\ \left| I\right| =p,\left| J\right| =q}}^{r}(-\sum_{i\in
(I\cap J)}\kappa _{i}+\sum_{i\in (I\cup J)^{c}}\kappa _{i})\left| \varphi _{I%
\overline{J}}^{a}\right| ^{2}
\end{equation*}
the remainder of the proof goes through as in \cite{kobayashi} .
\end{remark}

We now examine the consequences of the results in this section for the Higgs
bundles defined by \ref{higgsform}. With $(E,\theta )$ as in \ref{higgsform}
and assuming $(X,g)$ is a K\"{a}hler manifold, we get a second grading of
the bundle $E\bigotimes \Lambda ^{w}(X)=\bigoplus_{p+q=w}E\bigotimes \Lambda
^{(p,q)}(X)=\bigoplus_{\substack{ p+q=w  \\ 1\leq a\leq n}}E_{a}^{(p,q)}$
where $E_{a}^{(p,q)}=\Lambda ^{(a,0)}(X)\bigotimes \Lambda ^{(p,q)}(X)$.
Then one checks $\theta (\overline{\theta }_{h}^{\ast
}):E_{a}^{(p,q)}\rightarrow E_{a+k-1}^{(p+1,q)}(E_{a+k-1}^{(p,q-1)})$, and
if $E$ is endowed with a natural metric $h$ (\ref{hnatural}) then $\overline{%
\theta }_{h}(\theta ^{\ast }):E_{a}^{(p,q)}\rightarrow
E_{a-k+1}^{(p,q+1)}(E_{a-k+1}^{(p-1,q)})$.

We continue assuming $E$ is endowed with a natural Hermitian metric $h$.\
Then the decomposition $E=\bigoplus_{a=1}^{n}\Lambda ^{(a,0)}(X)$ is $h$
orthogonal, hence is also preserved by the associated Hermitian connection
and its curvature, $\Theta .$ It follows that $e(\Theta
):E_{a}^{(p,q)}\rightarrow E_{a}^{(p+1,q+1)}$, that $i(\Lambda e(\Theta
)-e(\Theta )\Lambda ):E_{a}^{(p,q)}\rightarrow E_{a}^{(p,q)},$ and that $%
i(\Lambda e(F_{h}^{(1,1)})-e(F_{h}^{\left( 1,1\right) })\Lambda )=i(\Lambda
e(\Theta )-e(\Theta )\Lambda )+\overline{\theta }_{h}\overline{\theta }%
_{h}^{\ast }+\overline{\theta }_{h}^{\ast }\overline{\theta }_{h}-\theta
\theta ^{\ast }-\theta ^{\ast }\theta :E_{a}^{(p,q)}\rightarrow
E_{a}^{(p,q)} $ (cf. \ref{kodnakcrv}). Finally if we assume that the Higgs
form $\theta $ is $h$ parallel then $F_{h}=F_{h}^{\left( 1,1\right) }\,$%
(Proposition \pageref{prop}) and we have $i(\Lambda e(F_{h})-e(F_{h})\Lambda
),\overline{\boxminus }=\square _{D^{\prime \prime }},\boxminus =\square
_{D_{h}^{^{\prime }}}:E_{a}^{(p,q)}\rightarrow E_{a}^{(p,q)}$. From this
data we deduce the following interpretation of the last part of Theorem \ref
{thm3.2} in this setting.

\begin{theorem}
\label{thmdoublehiggs}Let $(X,g)$ be a compact K\"{a}hler manifold of
complex dimension n, let $(E,\theta )$ $\rightarrow X$ be the Higgs bundle
given by \ref{higgsform}. Let $h$ be a natural metric on $E$ and assume that
the Higgs form $\theta $ is parallel, $d^{\bigtriangledown }\theta =0.$ For
a section $s$ of $E\bigotimes \Lambda ^{w}(X)$, write $s=\sum_{\substack{ %
p+q=w  \\ 1\leq a\leq n}}s_{a}^{(p,q)},$ where $s_{a}^{(p,q)}\in C^{\infty
}\Gamma E_{a}^{(p,q)}$. If $i(\Lambda e(F_{h})-e(F_{h})\Lambda )\leq 0$
pointwise as an operator on $E\bigotimes \Lambda ^{w}(X),$ then $\square
_{D^{\prime \prime }}s=0$ implies $0=$ $\square _{D^{\prime \prime
}}s_{a}^{(p,q)}=\square _{D_{h}^{^{\prime }}}s_{a}^{(p,q)}=\square
_{D_{h}}s_{a}^{(p,q)}$ for all $(p,q)$ and all $a$.
\end{theorem}

Using Proposition \ref{wparallel} we deduce the following

\begin{corollary}
Let $(X,g)$ be a compact K\"{a}hler manifold of complex dimension n, let $%
(E,\theta )$ $\rightarrow X$ be the Higgs bundle given by \ref{higgsform}.
Assume that $\varpi $ is $g$ parallel and use the standard extension of $g$
as the Hermitian metric on $E.$ For a section $s$ of $E\bigotimes \Lambda
^{w}(X)$, write $s=\sum_{\substack{ p+q=w  \\ 1\leq a\leq n}}s_{a}^{(p,q)},$
where $s_{a}^{(p,q)}\in C^{\infty }\Gamma E_{a}^{(p,q)}$. If $i(\Lambda
e(F_{h})-e(F_{h})\Lambda )\leq 0$ pointwise as an operator on $E\bigotimes
\Lambda ^{w}(X),$ then $\square _{D^{\prime \prime }}s=0$ implies $0=$ $%
\square _{D^{\prime \prime }}s_{a}^{(p,q)}=\square _{D_{h}^{^{\prime
}}}s_{a}^{(p,q)}=\square _{D_{h}}s_{a}^{(p,q)}$ for all $(p,q)$ and all $a$.
\end{corollary}

\end{document}